\title{\vspace{-1cm} Minors in expanding graphs}
\author{Michael Krivelevich\thanks{ Department of Mathematics,
Raymond and Beverly Sackler Faculty of Exact Sciences, Tel Aviv
University, Tel Aviv 69978, Israel. E-mail:
krivelev@post.tau.ac.il. Research supported in part by USA-Israel
BSF Grant 2002-133 and by grant 526/05 from the Israel Science
Foundation.} \and Benjamin Sudakov\thanks{
Department of Mathematics, Princeton University, Princeton, NJ 08544, and
Institute for Advanced Study, Princeton. E-mail:
bsudakov@math.princeton.edu.
Research supported in part by NSF CAREER award DMS-0546523, NSF grant
DMS-0355497, USA-Israeli BSF grant, Alfred P. Sloan fellowship, and
the State of New Jersey.
 }}
\date{}
\documentclass [11pt]{article}
\usepackage{graphics}
\usepackage{amsfonts}
\usepackage{amsmath}
\newtheorem{theorem}{Theorem}
\newtheorem{defin}{Definition}
\newtheorem{theo}{Theorem}[section]
\newtheorem{prop}[theo]{Proposition}
\newtheorem{propos}{Proposition}
\newtheorem{lemma}[theo]{Lemma}
\newtheorem{coro}[theo]{Corollary}

\renewcommand{\baselinestretch}{1.06}
\begin{document}
\maketitle

{\small \renewcommand{\baselinestretch}{0.7} \tableofcontents}

\section{Brief summary of results}\label{summary}
In this paper we address several extremal problems related to
graph minors. In all of our results we  assume essentially that a
given graph $G$ is expanding, where expansion is either postulated
directly, or $G$ can be shown to contain a large expanding
subgraph, or $G$ is locally expanding due to the fact that $G$
does not contain a copy of a fixed bipartite graph $H$. We need the
following definitions to state our results. A graph $\Gamma=(U,F)$ with
vertex set $U=\{u_1,\ldots,u_k\}$ is a {\em minor} of a graph
$G=(V,E)$ if the vertex set $V$ of $G$ contains a sequence of
disjoint subsets $A_1,\ldots, A_k$ such that the induced subgraphs
$G[A_i]$ are connected, and there is an edge of $G$ between $A_i$
and $A_j$ whenever the corresponding vertices $u_i,u_j$ of $\Gamma$ are
connected by an edge. A graph $G=(V,E)$ is {\em
$(t,\alpha)$-expanding} if every subset $X\subset V$ of size
$|X|\le \alpha |V|/t$ has at least $t|X|$ external neighbors in
$G$. A graph $G=(V,E)$ is called {\em $(p,\beta)$-jumbled} if
$$
\left| e(X)-p\frac{|X|^2}{2}\right|\le \beta |X|\
$$
for every subset $X\subseteq V$, where $e(X)$ stands for the
number of edges spanned by $X$ in $G$. Informally, this definition
indicates that the edge distribution of $G$ is similar to that of
the random graph $G_{|V|,p}$, where the degree of similarity is
controlled by parameter $\beta$.

Here are the main results of this paper.

\begin{theorem}\label{th1}
Let $0<\alpha<1$ be a constant. Let $G$ be a
$(t,\alpha)$-expanding graph of order $n$, and let $t\geq 10$.
Then $G$ contains a minor with average degree at least
$$
c \frac{\sqrt{nt \log t}}{\sqrt{\log n}},
$$
where $c=c(\alpha)>0$ is a constant.
\end{theorem}
This is an extension of results of  Alon, Seymour and Thomas \cite{AST}, Plotkin, Rao and Smith
\cite{PRS}, and of Kleinberg and Rubinfeld \cite{KR}, who cover basically the case
of expansion by a constant factor $t=\Theta(1)$.

\begin{theorem}\label{th2}
Let $G$ be a $(p,\beta)$-jumbled graph of order $n$ such that
$\beta=o(np)$. Then $G$ contains a minor with average degree
$cn\sqrt{p}$, for an absolute constant $c>0$.
\end{theorem}
This statement is an extension of results of A. Thomason \cite{T1, T2}, who
studied the case of constant $p$.  It can be also used 
to derive some of the results of Drier and Linial \cite{DL}.

\begin{theorem}\label{th3}
Let $2\le s\le s'$ be integers. Let $G$ be a $K_{s,s'}$-free graph
with average degree $r$. Then $G$ contains a minor with average
degree $cr^{1+\frac{1}{2(s-1)}}$, where $c=c(s,s')>0$ is a
constant.
\end{theorem}
This confirms a conjecture of K\"uhn and Osthus from \cite{KO}.

\begin{theorem}\label{th4}
Let $k \geq 2$ and let $G$ be a $C_{2k}$-free graph with average
degree $r$. Then $G$ contains a minor with average degree
$cr^{\frac{k+1}{2}}$, where $c=c(k)>0$ is a constant.
\end{theorem}
This theorem generalizes results of Thomassen
\cite{Th}, Diestel and Rompel \cite{DR}, and K\"uhn and Osthus \cite{KO1}, who
proved similar statements under the (much more restrictive) assumption
that $G$ has girth at least $2k+1$.

All of the above results are, up to a constant factor,
asymptotically tight (Theorems \ref{th1}, \ref{th2}), or are
allegedly tight (Theorems \ref{th3}, \ref{th4}), where in the
latter case the tightness hinges upon widely accepted conjectures
from Extremal Graph Theory about the asymptotic behavior of the
Tur\'an numbers of $K_{s,s'}$ and of $C_{2k}$.

\section{Background}\label{s-backg}
This paper is devoted to two of the most fundamental, yet normally
quite distant, concepts in modern Graph Theory -- minors and
expanding graphs. Their prominent role in mathematics is reflected
by the fact that both have been featured in a popular column
``What is...?" of the AMS Notices \cite{Moh}, \cite{Sar}. The
purpose of this section is to provide a basic information for both
of these concepts, and also for several related notions in Graph
Theory, relevant for this paper. Before going into technicalities, we 
would like to state notational agreements to be 
used in this paper. All graphs considered here are finite, without loops and without
multiple edges, unless stated explicitly otherwise.
Most of our notation  is rather standard and can be found in any textbook in Graph Theory. 
Here we define several less common pieces of notation, used
throughout the paper. 

Let $G=(V,E)$ be a graph. For a subset
$X\subseteq V$ we denote by $e_G(X)$ or simply by $e(X)$ the
number of edges of $G$ spanned by $X$, and by $N(X)$ the external
neighborhood of $X$:
$$
N(X) := \{u: u\not\in X, u\mbox{ has a neighbor in } X\}\ .
$$
In case $X=\{v\}$ we simply write $N(\{v\})=N(v)$; obviously, the
cardinality of $N(v)$ is the degree of $v$ in $G$. For two
disjoint sets $X,Y\subset V$, we denote the number of edges of $G$
connecting $X$ and $Y$ by $e(X,Y)$.

As quite customary in Extremal Graph Theory, our approach to the
problems researched will be asymptotic in nature. We thus assume
that an underlying parameter (normally the order $n$ of a graph)
tends to infinity and is therefore assumed to be sufficiently large
whenever necessary. We also do not make any serious attempt to optimize
absolute constants in our statements and proofs. 
All logarithms are in the natural basis. We omit systematically
rounding signs for the sake of clarity of presentation.

The following (standard) asymptotic notation
will be utilized extensively: for two functions $f(n)$, $g(n)$ of
a natural valued parameter $n$, we write $f(n)=o(g(n))$, whenever
$\lim_{n \rightarrow\infty} f(n)/g(n)=0$; $f(n)=O(g(n))$ if there
exists a constant $C>0$ such that $f(n)\le C g(n)$ for all $n$.
Also, $f(n)=\Omega(g(n))$ if $g(n)=O(f(n))$, and $f(n)=\Theta(g(n))$
if both $f(n)=O(g(n))$ and $f(n)=\Omega(g(n))$ are satisfied.

\subsection{Minors}
\begin{defin}\label{minor-def}
A graph $\Gamma$ is  a {\em minor} of a graph $G$ is for every vertex
$u\in \Gamma$ there is a connected subgraph $G_u$ of $G$ such that all
subgraphs $G_u$ are vertex disjoint, and $G$ contains an edge between
$G_u$ and $G_{u'}$ whenever $(u,u')$ is an edge of $\Gamma$.
\end{defin}
 An equivalent definition is through edge deletions and contractions:
we can obtain a minor $\Gamma$ of a graph $G$ by first deleting all
edges except those in subgraphs $G_u$, $u\in \Gamma$, and those
connecting $G_u$, $G_{u'}$ for $(u,u')\in E(\Gamma)$, and then
contracting all edges inside each of the connected subgraphs
$G_u$. (Given an edge $e=(v',v'')$ of a graph $G$, {\em contracting}
$e$ results in replacing $v',v''$ by a single new vertex $v$, and
connecting $w\in V(G)-\{v',v''\}$ to the new vertex $v$ if and only
if $w$ is connected to $v'$ or to $v''$ or to both in $G$).

Though the notion of graph minors appears at the first sight to be
purely graph theoretic, it turns out to be absolutely essential in
bridging between Graph Theory on one side, and Topology and
Geometry on the other -- one of the most non-trivial and
fundamental connections in Mathematics. Indeed, the famous theorem
of Kuratowski \cite{Kur30} (in its reformulation due to Wagner
\cite{Wag37}) postulates that a graph $G$ can be embedded in the
plane (is planar) if and only if neither the complete graph $K_5$
on five vertices nor a complete bipartite graph $K_{3,3}$ with
three vertices at each side is a minor of $G$. This was the
beginning of Topological Graph Theory, whose culmination is
without a doubt the celebrated Robertson-Seymour theory of graph
minors. In a series of twenty papers, spanning over two decades
(with \cite{RS-XX} being the concluding paper of the series),
Robertson and Seymour proved the so called Wagner conjecture: in
every infinite collection of graphs, there are two such that one
is a minor of the other (in other words, the set of finite graphs
with the ``minor" relation as partial order is well-quasi-ordered;
see, e.g. \cite{Kru72} for more information about the theory of
well-quasi-ordering). An equivalent formulation is that every
family of graphs closed with respect to taking minors can be
characterized by a finite family of excluded minors. As a
corollary Robertson and Seymour were able to derive that for every
closed compact surface there is a finite list of graphs such that
a graph $G$ is embeddable in this surface if and only if it does
not contain any of these as a minor. This is of course an
extremely far-reaching generalization of the Kuratowski theorem.
The Robertson-Seymour Structural Graph Theory is undoubtedly an
admirable research effort and one of the crown achievements of
Combinatorics, whose impact is truly immense. As our research in
minors will proceed along rather different lines, we will not
dwell on this wonderful theory anymore, referring the reader
instead to a very nice survey of Lov\'asz on the subject
\cite{Lov05}.

\subsection{Expanding graphs}
The second fundamental concept of this paper is expanding graphs.
Informally, a graph $G$ is said to be an expanding graph or an
expander if every subset $X$ of $V(G)$ has relatively many
neighbors outside $X$. (This is what is usually called {\em vertex
expansion}, sometimes an alternative notion of {\em edge
expansion} is used, there every set $X$ is required to be incident
to many edges crossing between $X$ and its complement in $V(G)$;
for constant degree graphs these two notions are essentially
equivalent). Of course, a formal definition is required here,
firstly, to measure the expansion quantitatively, and secondly to
distinguish between the expansion of small and large sets --
obviously a set $X$ containing half the vertices of $V$ cannot
have more than $|X|$ outside neighbors, while a much smaller set
$X$ can expand by a much larger factor. There are several
definitions of expanders in common use, capturing sometimes rather
different expansion properties. In this paper we find it much more
important to look at the expansion of small sets, and for this
reason we adopt the following formal definition of an expander.

\begin{defin}\label{def-exp}
Let $t>0$, $0<\alpha <1$. A graph $G=(V,E)$ is {\em
$(t,\alpha)$-expanding} if every subset $X\subset V$ of size
$|X|\le \alpha |V|/t$ has at least $t|X|$ external neighbors in
$G$.
\end{defin}
Normally we will think of $\alpha$ as being an absolute constant.
In this case, the above definition says that every set $X$ of size
$|X|=O(n/t)$ expands by a factor of at least $t$.

As the research in the last quarter century has convincingly
shown, the notion of expanders is of utmost value in an amazing
variety of fields, both in and outside of Discrete Mathematics.
Applications include design of efficient communication networks,
error-correcting codes with efficient encoding and decoding,
derandomization of randomized algorithms, study of metric embeddings,
to mention just a few.
Expanders are usually constructed much easier using probabilistic,
existential arguments (see, e.g. \cite{Pinsker}); explicit
constructions of expander graphs are much harder to come by and
range from classical papers of Margulis \cite{Mar88} and of Lubotzky,
Phillips and Sarnak \cite{LPS88}, to a relatively recent zig-zag
product construction of Reingold, Vadhan and Wigderson
\cite{RVW02}.

Our viewpoint here will be somewhat different from the above
mentioned papers. Instead of discussing ways to construct good
expanders, we will concentrate on {\em properties} of expanders,
and more specifically on the appearance of large minors in
expanding graphs.

General information about expanders, their properties and applications can be
found in a recent excellent survey of Hoory, Linial and 
Wigderson \cite{HLW}.

\subsection{Pseudo-random graphs}
A notion closely related to expanding graphs is that of
pseudo-random graphs. As the name clearly suggests, pseudo-random
graphs can be informally described as graphs resembling truly
random graphs, most commonly the so called binomial random graphs
$G_{n,p}$. We first remind the reader the definition of this
probability space. Given two parameters $n$ and $0\le p\le 1$, the
{\em random graph} $G_{n,p}$ is a probability space of all graphs
on $n$ vertices labeled $1,...,n$, where for each pair $1\le i\ne
j\le n$, the probability that $(i,j)$ is an edge is $p$,
independently of all other pairs. Equivalently, $G_{n,p}$ is the
probability spaces of all labeled graphs with vertex set
$\{1,\ldots,n\}$, endowed with the probability measure
$Pr[G]=p^{|E(G)|}(1-p)^{{n\choose 2}-|E(G)|}$. In quite a few
cases the edge probability $p$ is in fact a function $p=p(n)$ of
the number of vertices $n$, vanishing as $n$ tends to infinity.
We say that random graph possesses a property $\cal P$ {\em with 
high probability }, if the probability that $G_{n,p}$ satisfies $\cal P$ 
tends to 1 as $n$ tends to infinity. 
This probability space is undoubtedly the most studied and the most convenient to
work with probability distribution on graphs. When defining
pseudo-random graphs, one usually tries to capture quantitatively
their similarity to truly random graphs, in this aspect or
another. Arguably the most important feature of random graphs is
their edge distribution, and so it is quite natural to expect that
a definition of a pseudo-random graph will address this property.
For the probability space $G_{n,p}$, edge distribution is quite
easy to handle -- for a given subset $X\subseteq V(G)$, the number
of edges spanned by $X$ in $G_{n,p}$ is a binomially distributed
random variable with parameters ${{|X|}\choose {2}}$ and $p$;
applying standard bounds on the tails of the binomial distribution
one can easily show that with high probability {\em all} sets $X$
of cardinality $k$ span indeed close to ${k\choose 2}p$ edges in
$G_{n,p}$, if $k$ is not too small. This fact
motivates the following definition of a pseudo-random graphs due
to Thomason \cite{T1}, \cite{T2}:

\begin{defin}\label{def-jumbled}
A graph $G=(V,E)$ is {\em $(p,\beta)$-jumbled} if for every subset $X\subseteq V(G)$,
$$
\left|e_G(X)-\frac{p|X|^2}{2}\right|\le \beta |X|\,.
$$
\end{defin}

Thus, if $G$ is a $(p,\beta)$-jumbled graph, its edge density is
around $p$, and its edge distribution is similar to that of the random
graph $G_{n,p}$, where the degree of similarity (or rather of
proximity to the expected number of edges) is controlled by the
parameter $\beta$. Random graphs $G_{n,p}$ are easily shown to be
$(p,O(\sqrt{np}))$-jumbled for all not too small values of the edge
probability $p$. Moreover, one can show (see \cite{ES}) that if a
graph $G$ on $n$ vertices is $(p,\beta)$-jumbled, then
$\beta=\Omega(\sqrt{np})$; for this reason $(p,\beta)$-jumbled
graphs $G$ with $\beta=\Theta(\sqrt{np})$ are considered very good
pseudo-random graphs.

Pseudo-random graphs is a central concept in modern Combinatorics,
whose importance is derived in part from that of random graphs.
Quite a few known constructions of pseudo-random graphs are
deterministic, allowing thus to substitute somewhat elusive truly
random graphs, defined through probabilistic, existential means,
with quite accessible deterministic descriptions -- a feature
crucial in a variety of applications. Moreover, in certain
applications one can utilize features of (carefully crafted)
pseudo-random graphs non-existent typically in random graphs of
the same edge density.

As we have indicated, several alternative definitions of
pseudo-random graphs are available; here we describe just one of
them, based on graph spectrum. Given a graph $G=(V,E)$ with vertex
set $V=\{v_1,\ldots,v_n\}$, the {\em adjacency matrix} of $G$ is
an $n$-by-$n$ matrix $A$ of zeroes and ones, defined by:
$a_{ij}=1$ if and only if $(v_i,v_j)\in E(G)$, and $a_{ij}=0$
otherwise. It is easy to observe that $A$ is a symmetric real
matrix, and therefore $A$ has a full set of $n$ real eigenvalues,
denoted by $\lambda_1,\ldots,\lambda_n$, customarily sorted in the
non-increasing order $\lambda_1\ge\lambda_2\ge\ldots\ge \lambda_n$
and usually called the {\em eigenvalues} of the graph $G$ itself.
If $G$ is a $d$-regular graph, then the first eigenvalue
$\lambda_1$ is easily seen to be $\lambda_1=d$ (with the
corresponding eigenvector being the all-one vector), while all
others satisfy $|\lambda_i|\le d$, $i=2,\ldots,n$. Now, equipped
with this terminology, we can give an alternative definition of a
pseudo-random graph introduced by Alon. A graph $G=(V,E)$ is called an {\em
$(n,d,\lambda)$-graph} if $G$ has $n$ vertices, is $d$-regular,
and in addition all of its eigenvalues but the first one satisfy:
$|\lambda_i|\le \lambda$, $i=2,\ldots,n$. A very frequently used
result from Spectral Graph Theory (see, e.g., Chapter 9 in \cite{AS}) 
postulates that if $G$ is an $(n,d,\lambda)$-graph, then
$$
\left| e_G(X)-\frac{p|X|^2}{2}\right| \le \lambda |X|\,,
$$
for all subsets $X\subseteq V(G)$, implying that an
$(n,d,\lambda)$-graph is $(d/n,\lambda)$-jumbled. Several
constructions of $(n,d,\lambda)$-graphs with $\lambda=O(\sqrt{d})$
are available, they are based on a variety of algebraic and
geometric properties. We would like to mention in passing that
graph eigenvalues are frequently used to ensure graph expansion
too.

The reader is advised to consult a survey \cite{KS} on
pseudo-random graphs by the authors for an extensive coverage of
pseudo-random graphs, their definitions and properties.

\section{Extremal problems for minors}\label{s-extrem}
The subject of this paper can be classified as ``Extremal problems
for minors". Given the prominence of these two branches of Graph
Theory (theory of minors and extremal graph theory), it is quite
natural to expect the appearance of results combining these two
subjects. And indeed, our paper is certainly not the first to
address extremal problems for minors; in fact, this is already a
well established part of Graph Theory, with a variety of results
achieved. A recent survey of Thomason \cite{Th-sur} on the subject
describes several of its achievements.

Generally speaking, the motto of the extremal minor theory can be
stated as finding sufficient conditions for the existence of a
minor from given family, or a concrete minor (say, a clique minor
or certain order) in a given graph. Here is an illustrative
example of a result of this sort: every graph $G$ on $n$ vertices
with more than $3n-6$ edges contains a complete graph $K_5$ or a
complete bipartite graph $K_{3,3}$ as a minor. This is of course
nothing else but rephrasing of the Kuratowski-Wagner theorem
combined with the classical fact that a planar graph on $n$
vertices has at most $3n-6$ edges (which in turn follows easily
from the celebrated Euler formula connecting the numbers of
vertices, edges and faces in any planar embedding). As yet another
illustration we can mention the famous Hadwiger Conjecture,
suggesting that a graph that cannot be properly colored with $k$
colors has a clique $K_{k+1}$ as a minor; this notorious
conjecture has been proven so far for very few initial values of
$k$, see \cite{Toft} for a survey of its status.

Here we will be mostly looking for results of the following sort:
if a graph $G$ is sufficiently dense, or has sufficiently large
average degree (plus possibly additional conditions imposed),
then $G$ contains a large minor. Perhaps the best known result of
this sort was proved independently by Kostochka \cite{K} and by
Thomason \cite{T3}, who showed that there exists an absolute
constant $c>0$ such that every graph $G$ with average degree
$d=2|E(G)|/|V(G)|$ contains a clique on $cd/\sqrt{\log d}$
vertices as a minor. Recently the asymptotic value of $c$ has been
determined by Thomason \cite{T4}.

Under certain additional conditions one can guarantee a clique
minor of order (much) larger than $d/\sqrt{\log d}$ in a graph of
average degree $d$. Several of our results are indeed of this
type. When looking for large minors one should remember however
that there is a limit of the size of a minor one can find in a
graph. This limit is given by the following very simple yet very
useful observation.

\begin{propos}\label{propos1}
Let $H$ be a minor of $G$. Then the number of edges of $H$ does
not exceed the number of edges of $G$.
\end{propos}

The above proposition immediately implies that a graph $G$ on $n$
vertices with average degree $d$ cannot contain a graph $\Gamma$ with 
average degree $k>\sqrt{nd}$ as minor. Indeed, the number of edges of
$G$ is $nd/2$, and thus if $\Gamma$ is a minor of $G$ then $\frac{k^2}{2}\le \frac{nd}{2}$. 
We will repeatedly use this simple bound as a
benchmark to measure the quality of our results.

In the rest of this section we survey a variety of known results
in Extremal Minor Theory, having in mind our results and their
comparison to the previously obtained results.

There are several results connecting between (the absence of)
separators and minors in graphs. A {\em separator} $S$ of a graph $G$ is
a set of vertices whose removal separates the graph into connected
components, each of size at most $\frac{2}{3}|V(G)|$. Alon,
Seymour and Thomas \cite{AST} proved that a graph of order $n$
without a $K_h$ minor has a separator of size $O(h^{3/2}n^{1/2})$.
This was extended to large $h$ by Plotkin, Rao and Smith \cite{PRS} who proved
that a graph without a $K_h$ minor has a separator of size $O(h
\sqrt{n \log n})$. The last result implies in particular that an
expander graph of constant degree has a clique minor of size
$\Omega(\sqrt{n/\log n})$. On the other hand, since every graph
has trivially a separator of size $n/3$, one can only show the
existence of a clique minor of order at most
$O(\sqrt{n/\log n})$ using these results.

Kleinberg and Rubinfeld addressed in \cite{KR} a connection
between expansion and the existence of large minors. They used the
following, rather weak, definition of expansion: a graph $G$ is an
$\alpha$-expander if every set $X$ of at most half of the vertices
of $G$ has at least  $\alpha|X|$ outside neighbors in $G$. It is
proven in  \cite{KR} that for every fixed $\alpha>0$ there is a
constant $c>0$ such that an $\alpha$-expander graph of order $n$
contains every graph $H$ with at most $n/\log^c n$ vertices and edges
as a minor. While this result is quite useful in finding
large minors in sparse graphs (in particular those of constant
maximum degree), it appears to be of rather limited value for the
denser case and can not be used to show the existence of
a clique minor of order larger than $\Omega(\sqrt{n/\log n})$.

Sunil Chandran and Subramanian \cite{CS} discussed a connection
between spectral properties of a graph and its minors. They proved
in particular that if $G$ is a $d$-regular graph on $n$ vertices whose
second eigenvalue is at most $\lambda$, then $G$ contains a clique
minor on
$\Omega\left(\big(\frac{n(d-\lambda)^2}{(3d-2\lambda)^2}\big)^{1/3}\right)$
vertices. Observe that this result can be used only to show the
existence of clique minors of order up to $cn^{1/3}$, which is a
relatively weak bound.

Another avenue of research in extremal problems in minors (also
pursued in this paper) aims to prove the existence of large minors
in graphs with excluded subgraphs. K\"uhn and Otshus proved in
\cite{KO} that for all integers $2\le s\le s'$ there exist
constants $r_0=r_0(s,s')$ and $c=c(s,s')$ such that every
$K_{s,s'}$-free graph $G$ of average degree $r\ge r_0$ contains a
minor of average degree $d$ satisfying
$$
d \ge c\frac{r^{1+\frac{1}{2(s-1)}}}{(\log r)^{2+\frac{1}{s+1}}}
\,.
$$
They conjectured
however that the logarithmic term is not needed in this bound and were
able to verify this conjecture for the case when the graph $G$ is
assumed to be regular. Observe that after having obtained a minor
of average degree $d$ one can use the above mentioned results of
Kostochka and Thomason \cite{K}, \cite{T3}, \cite{T4} to derive
the existence of a clique minor on $cd/\sqrt{\log d}$ vertices.

Another nice result of K\"uhn and Osthus guarantees the
existence of large minors in graphs with large girth (i.e. without
short cycles). They proved in \cite{KO1} that for every odd
integer $g\ge 5$ there exists a constant $c=c(g)>0$ such that
every graph $G$ of average degree $r$ and without cycles shorter
than $g$ (such a graph is said to have {\em girth} more than $g$)
contains a minor with average degree at least $cr^{(g+1)/4}$. This
result improves significantly a much earlier result of Thomassen
\cite{Th} and a recently obtained result by Diestel
and Rompel \cite{DR}. Observe that the assumption for the case
$g=5$ essentially amounts to forbidding a 4-cycle, or $K_{2,2}$;
thus this result of K\"uhn and Osthus establishes their above
mentioned conjecture for the case $s=s'=2$.

Bollob\'as, Catlin and Erd\H{o}s \cite{BCE} analyzed the
appearance of large minors in random graphs. They proved that for
a constant edge probability $p$, $0<p<1$, the largest clique minor
in a random graph $G_{n,p}$ is of order $n/\sqrt{\log n}$ (in fact,
their result is more accurate -- they were able to establish not
only the asymptotic order of magnitude of the largest clique minor
in $G_{n,p}$, but actually its asymptotic behavior). As a result,
and taking into account a well known fact that the chromatic
number of $G_{n,p}$ in this range is with high probability
$O(n/\log n)$, Bollob\'as et al. were able to derive that almost
every graph satisfies the Hadwiger conjecture. The argument of
\cite{BCE} can be used to show that with high probability the
largest clique minor in $G_{n,p}$ has order of magnitude
$\Theta\big(n\sqrt{p}/\sqrt{\log n}\big)$, for subconstant values of the edge
probability $p(n)$ as well.

Much less is known in the case of pseudo-random (or jumbled) graphs.
Thomason proved in \cite{T1} (see also \cite{T2}) that $(p, \beta)$-jumbled graphs
with $p$ constant and $\beta=O(n^{1-\epsilon})$ contain a clique
minor of size at least $(1+o(1))n/\sqrt{\log_b n}$, where $b=1/(1-p)$.
For small $p$, this has the same order of magnitude $\frac{n\sqrt{p}}{\sqrt{\log n}}$ 
as the result for $G_{n,p}$.

Finally, we mention a recent result of Drier and Linial \cite{DL}
who discussed minors in lifts of graphs. An $\ell$-lift of a
labeled graph $G=(V,E)$ is a graph with vertex set
$V\times[\ell]$, whose edge set is the union of perfect matchings
between $\{u\}\times[\ell]$ and $\{v\}\times[\ell]$ for each edge
$(u,v)\in E$. In a random lift these matchings are selected
uniformly at random. Drier and Linial proved that for $\ell \leq
O(\log n)$ almost every lift of the complete graph $K_n$ contains
a clique minor of size $\Theta(n)$, and for $\ell >\log n$ it
contains a clique minor of size at least
$\Omega\Big(\frac{n\sqrt{\ell}}{\sqrt{\log (n\ell)}}\Big)$. The last result was
shown to be tight in \cite{DL} as long as $\log n<
\ell<n^{1/3-\epsilon}$.

\section{Our results}\label{s-results}
In this section we present in full details the results of this
paper. We also compare them with previously obtained results,
surveyed in brief in Section \ref{s-extrem}, and discuss their
tightness.

The first of our results is about minors in expanding graphs. We
prove:

\begin{theo}
\label{minor-t-expand} Let $G$ be a $(t, \alpha)$-expanding graph
of order $n$ and let $t\geq 10$. Then $G$ contains a minor with
average degree at least
$$
c\alpha^3 \frac{\sqrt{nt \log t}}{\sqrt{\log n}}\,,
$$
where $c>0$ is some absolute constant independent of $\alpha$.
\end{theo}

This theorem together with  the results of Kostochka \cite{K} and
Thomason \cite{T3} mentioned in Section \ref{s-extrem} gives the
following corollary.

\begin{coro}
\label{cliqueminor-t-expand} Let $G$ be a $(t, \alpha)$-expanding
graph of order $n$, and let $t\geq 10$. Then $G$ contains a clique
minor of size
$$
c\alpha^3 \frac{\sqrt{nt \log t}}{\log n}\,,
$$
where $c$ is some absolute constant independent of $\alpha$.
\end{coro}

For $t\geq n^{\epsilon}$ this gives a clique minor of size
$\Omega\left(\frac{\sqrt{nt}}{\sqrt{\log n}}\right)$. The random
graph $G_{n,p}$ with $p=10t/n$ can be easily shown with high probability to be 
a $(t,0.5)$-expander in this range of $t$, and as we mentioned
before its largest clique minor is typically of order
$O\left(\frac{\sqrt{nt}}{\sqrt{\log n}}\right)$. This shows that
our result is tight up to a constant factor. For small values of
$t\leq \log n$ the result of this corollary can be slightly
improved as follows:

\begin{prop}
\label{small-t} If $G$ is a $(t, \alpha)$-expanding graph of order
$n$ and $t\geq 10$, then $G$ contains a clique minor of size
$$
\Omega\left(\alpha^{2} \sqrt{\frac{n \log t}{\log n}}\right)\,.
$$
\end{prop}

Observe that the above results constitute a substantial extension
of the results of  Alon, Seymour and Thomas \cite{AST}, Plotkin, Rao and Smith \cite{PRS}, and of
Kleinberg and Rubinfeld \cite{KR}, that cover basically the case
of expansion by a constant factor $t=\Theta(1)$. Our results,
though applicable also for the case $t=\Theta(1)$, enable to show
the existence of larger minors whenever the expansion factor $t$
becomes super-constant.

The next our result is about minors in pseudo-random (or jumbled)
graphs. We prove:

\begin{theo}
\label{jumbled-minor} Let $G$ be a $(p,\beta)$-jumbled graph of
order $n$ such that $ \beta=o(np)$. Then $G$ contains a minor with
average degree $\Omega(n\sqrt{p})$.
\end{theo}

This statement is an extension of the results of Thomason \cite{T1, T2}, who
studied the case of constant $p$. As a $(p,\beta)$-jumbled graph $G$ on $n$ vertices with
$\beta=o(np)$ has average degree close to $np$ and thus
$\Theta(n^2p)$ edges, Proposition \ref{propos1} shows that Theorem
\ref{jumbled-minor} is asymptotically tight, up to a constant
factor. The above theorem also implies that an
$(n,d,\lambda)$-graph $G$ with $\lambda=o(d)$ has a minor of
average degree $\Omega(\sqrt{nd})$. This can be used in particular
to derive some of the results of Drier and Linial \cite{DL} on
minors in random lifts. For example, when $ \ell \ll n$ we have that
with high probability every pair of vertices in a random
$\ell$-lift of the complete graph $K_n$ has at most $(1+o(1))n/\ell$ common neighbors.
Using this one can easily show that this graph is an $(n\ell,n-1,\lambda)$-graph
with $\lambda=o(n)$. Therefore a random $\ell$-lift of  $K_n$ contains a minor with average degree
$\Omega(n\sqrt{\ell})$ and thus a clique minor of order $\Omega\Big(\frac{n\sqrt{\ell}}{\sqrt{\log (n\ell)}}\Big)$, 
by applying again Kostochka-Thomason. For larger values of $\ell$ one can obtain similar lower bound on the 
size of the clique minor in a random $\ell$-lift of the complete graph $K_n$ by first proving that 
all subsets of order at most $O(\ell)$ in such graph expand by a factor of $\Omega(n)$ and then using
Theorem \ref{minor-t-expand}.

The next group of results guarantees the existence of large minors
in graphs with excluded subgraphs. First, we prove:

\begin{theo}
\label{K_{s,s'}-minor} Let $2\le s\le s'$ be integers. Let $G$ be
a $K_{s,s'}$-free graph with average degree $r$. Then $G$ contains
a minor with average degree
$\Omega\left(r^{1+\frac{1}{2(s-1)}}\right)$.
\end{theo}

This confirms a conjecture of K\"uhn and Osthus from \cite{KO}.
The result is asymptotically tight modulo a well known and widely
accepted conjecture on the Tur\'an numbers of complete bipartite
graphs $K_{s,s"}$, saying that for constant $2\le s\le s'$, there
exists a $K_{s,s'}$-free graph $G$ on $n$ vertices with at least
$\Omega\big(n^{2-1/s}\big)$ edges. Denoting the average degree of such a
graph by $r$, we have then $r=\Omega(n^{1-1/s})$, and therefore by
Proposition \ref{propos1} a minor $H$ of $G$ has
$O\big(r^{2+1/(s-1)}\big)$ edges, and hence the average degree of $H$ is
at most $O\Big(r^{1+\frac{1}{2(s-1)}}\Big)$. The latter conjecture has been
settled for $s=2,3$ and all $s'\ge s$ (see, e.g., Chapter VI of
\cite{Bol}), furthermore, Alon, R\'onyai and Szab\'o proved it
\cite{ARS} for $s'> (s-1)!$; the asymptotic tightness of Theorem
\ref{K_{s,s'}-minor} thus follows in all these cases.

Theorem \ref{K_{s,s'}-minor} can be generalized somewhat to the
case where an excluded graph $H$ is a bipartite graph with bounded
degrees at one side. The corresponding result is:

\begin{theo}
\label{H-minor} Let $H$ be a bipartite graph of order $h$ with
parts $A$ and $B$ such that the degrees of all vertices in $B$ do
not exceed $s$. If $G$ is an $H$-free graph with average degree
$r$, then $G$ contains a minor with average degree
$\Omega\left(r^{1+\frac{1}{2(s-1)}}\right)$.
\end{theo}

Finally, we prove a minor-related result for $C_{2k}$-free graphs.

\begin{theo}
\label{C_2k-minor} Let $k \geq 2$ and let $G$ be a $C_{2k}$-free
graph with average degree $r$. Then $G$ contains a minor with
average degree $\Omega\left(r^{\frac{k+1}{2}}\right)$.
\end{theo}

This generalizes a result of K\"uhn and Osthus \cite{KO1}, who
proved such a theorem under the (much more restrictive) assumption
that $G$ has girth at least $2k+1$. Here too the asymptotic
optimality of Theorem \ref{C_2k-minor} relies on a well known
conjecture from Extremal Graph Theory (see, e.g., \cite{Bol}, p.
164), postulating that for any fixed $k\ge 2$, there exists a
graph $G$ on $n$ vertices without cycles of length up to $2k$ and
with $\Omega(n^{1+1/k})$ edges. This conjecture has been proven so
far for very few values of $k$.

Of course the Kostochka-Thomason result can be utilized to convert
minors with large average degree into clique minors, just as we
have done several times already.

The alert reader has probably noticed that in all three results
above the excluded fixed graph is bipartite. This is for a good
reason -- the complete bipartite graph $K_{r,r}$ is $H$-free for
any non-bipartite graph $H$, and yet every minor $H$ of $K_{r,r}$
has obviously average degree $O(r)$. This indicates that if one's
aim is to force an untypically large minor by excluding a fixed
graph $H$, $H$ should better be bipartite.

The rest of the paper is organized as follows. In Section
\ref{s-expanders} we discuss minors in expanding graphs and prove
Theorem \ref{minor-t-expand} and Proposition \ref{small-t}.
Section \ref{s-pseudo-random} is devoted to minors in
pseudo-random graphs, there we prove Theorem \ref{jumbled-minor}.
In Section \ref{s-H-free} we derive Theorems \ref{K_{s,s'}-minor}
and \ref{H-minor} about minors in $K_{s,s'}$-free graphs and in
$H$-free graphs. In Section \ref{s-C_2k-free} we prove Theorem
\ref{C_2k-minor} about large minors in $C_{2k}$-free graphs.
Section \ref{s-final} is devoted to concluding remarks.

\section{Minors in expanding graphs}\label{s-expanders}

In this section we prove Theorem \ref{minor-t-expand} and
Proposition \ref{small-t}.

Observe first that if $G$ is a $(t, \alpha)$-expanding graph of
order $n$, then every subset $X$ of $G$ of size $\alpha n/t\leq
|X| \leq \alpha n/2$ has $|N(X)|\geq \alpha n/2$. Indeed, such $X$
contains a subset $Y$ of size exactly $\alpha n/t$, hence
$|N(X)|\geq |N(Y)|-|X|\geq t|Y|-|X|\geq \alpha n/2$.

\begin{lemma} \label{diameter} Let $G$ be a connected $(s,
\beta)$-expanding graph of order $n$. Then the diameter of $G$ is
at most $3 \beta^{-1} \log n/\log s$.
\end{lemma}

\noindent {\bf Proof.} From the expansion of $G$ we have that for
every vertex $v$ and integer $q$ there are at least $\min\{s^q,
\beta n\}$ vertices which are within distance at most $q$ from
$v$. Taking $q=\log n/\log s$ we obtain that there are at least
$\beta n$ vertices within distance at most $\log n/\log s$ from
every vertex in $G$.

Now, suppose $G$ contains a pair of vertices $u,w$ such that the
distance between them is at least $3 \beta^{-1} \log n/\log s$.
Then on a shortest path from $u$ to $w$ we can find vertices
$v_1=u, \ldots, v_k=w$ such that $k>1/\beta$ and the distance
between every pair $v_i, v_j$ is at least $2\log n/\log s$. Denote
by $U_i$ the set of vertices which are at distance at most $\log
n/\log s$ from $v_i$. These sets are disjoint, each has size at
least $\beta n$ and therefore the size of their union is larger
than $n$. This contradiction completes the proof. \hfill $\Box$

\noindent
{\bf Proof of Theorem \ref{minor-t-expand}.}\, Let
$$p=\frac{\alpha^2}{100}\frac{\sqrt{nt \log t}}{\sqrt{\log n}} \quad \mbox{and}\quad
q=6\alpha^{-1}\frac{\sqrt{n\log n}}{\sqrt{t \log t}},$$ and
consider the following iterative procedure which we will repeat
$p$ times. In the beginning of iteration $k+1$ we will have $k$
disjoint sets $B_1, \ldots, B_k$ each of size $|B_i|=q$, such that
all induced subgraphs $G[B_i]$ are connected. We will construct a
new subset $B_{k+1}$, also of size $q$, such that induced subgraph
$G[B_{k+1}]$ is connected and there are at least $\alpha k/8$
indices $1 \leq i\leq k$ such that there is an edge from $B_i$ to
$B_{k+1}$. In the end of this algorithm if we contract all subsets
$B_i$ we will get a graph with average degree
$$\Omega(\alpha p)=\Omega\left(\alpha^3 \frac{\sqrt{nt \log t}}{\sqrt{\log n}}\right).$$

Let $B=\cup_{i=1}^k B_i$ and note that $|B|=b\leq pq \leq
0.06\alpha n$. Denote by $C=V(G)-B$ and by $G'$ the subgraph of
$G$ induced by $C$. Let $X$ be the subset of $C$ such that $2b/t
\leq |X|\leq \alpha n/t$ and $|N_{G'}(X)| < t|X|/2$. Then we have
$$|N_G(X)|\leq |N_{G'}(X)|+|B| \leq t|X|/2+b \leq t|X|, $$
which contradicts the assumption that $G$ is $(t,
\alpha)$-expanding. Therefore there exists $X \subset C $ of size
at most $2b/t$ such that the remaining set $D=C-X$ spans a
subgraph of $G$ in which every subset of size at most $\alpha n/t$
expands by a factor of at least $t/2$. Denote by $G''$ the
subgraph of $G$ induced by $D$. This graph might be disconnected,
but as we will see next it must have few very large components
which cover almost all its vertices.

Let $Y$ be a subset of  $G''$ such that  $3b/t\leq |Y|< \alpha
n/2$. Then, by the remark in the beginning of this section, we
have that $|N_G(Y)| \geq \min\{3b, \alpha n/2\} >|B|+|X|$. Hence
$Y$ has neighbors inside $D-Y$ and can not be an isolated
component of $G''$. Thus $G''$ contains a subset $Y$ of size at
most $3b/t$ such that all the vertices of $G''-Y$ are contained in
connected components of size at least $\alpha n/2$. Denote these
connected components by $G_1, \ldots, G_{\ell}$. Then clearly
$\ell \leq 2/\alpha$, and we also have that every subset of $G_i$
of size at most $\alpha n/t$ expands by factor at least $t/2$. By
Lemma \ref{diameter} (with $\beta=\alpha/2$ and $s=t/2$) this
implies that the diameter of each $G_i$ is at most $7
\alpha^{-1}\log n/\log t$.

Next we claim that there is an index $i$ such that there are at
least $r=\frac{k}{2\ell}$ sets $B_j$, each having at least
$\frac{t|B_j|}{2\ell}$ neighbors in $G_i$. If this is not the case
then we have $k-\ell \,\frac{k}{2\ell}=k/2$ sets $B_j$, each
having at most $\ell\,\frac{t|B_j|}{2\ell}= tq/2$ neighbors inside
$\cup_i G_i$. First suppose that $kq/2\leq \alpha n/t$. Then
taking a union of $k/2$ such sets $B_j$ we obtain a set of size
$b/2$ with at most $tb/4$ neighbors in $\cup_i G_i$. On the other
hand, by expansion the number of neighbors of this set in $G$ is
at least $tb/2$. Therefore the remaining $tb/4$ neighbors should
be inside $X\cup Y \cup B$. But this set has size at most $b+5b/t<
tb/4$, a contradiction. If $kq/2\geq \alpha n/t$ then we can take
a union of $\alpha n/(tq)$ such sets $B_j$ and obtain a set of
size $\alpha n/t$ with at most $(tq/2)(\alpha n/(tq))=\alpha n/2$
neighbors in $\cup_i G_i$. Again, by expansion, this set has at
least $\alpha n$ neighbors in $G$, so at least $\alpha n/2$ of
them are in $X\cup Y \cup B$. But the size of this set is not big
enough, a contradiction.

Therefore, without loss of generality, we can assume that each of
the first $r=\frac{k}{2\ell}$ sets $B_1, \ldots, B_r$ has at least
$t|B_j|/(2\ell)=tq/(2\ell)$ neighbors in $G_1$. Denote these sets
of neighbors by $U_1, \ldots, U_r$ respectively. Pick uniformly at
random with repetition $\frac{|G_1|}{tq/(2\ell)}$ vertices of
$G_1$ and denote this set by $W$. For every index $1 \leq i \leq
r$, the probability that $W$ does not intersect $U_i$ is at most
$\left(1-\frac{|U_i|}{|G_1|}\right)^{|W|}\leq 1/e$. Therefore the
expected number of sets $U_i$ which have non-empty intersection
with $W$ is at least $(1-1/e)r>r/2$. Hence there is a choice of
$W$ that intersects at least $r/2\geq k/(4\ell)\geq \alpha k/8$
sets $U_i$. Fix an arbitrary vertex  $w_0 \in W$ and consider a
collection of shortest paths in $G_1$ from $w_0$ to the remaining
vertices in $W$. Since the diameter of $G_1$ is at most $7
\alpha^{-1}\log n/\log t$ and
$$7 \alpha^{-1} |W|\log n/\log t \leq 7 \alpha^{-1} \frac{n}{tq/(2\ell)}\frac{\log n}{\log t}
\leq \frac{14}{3} \alpha^{-1} \frac{\sqrt{n\log n}}{\sqrt{t \log
t}}<q,$$ by taking a union of these paths  and adding extra
vertices if necessary we can construct a connected subset of size
$q$ containing $W$. Denote this set by $B_{k+1}$ and note that it
is connected by an edge to at least $\alpha k/8$ sets $U_i$,
$1\leq i \leq k$. This completes the proof of the theorem. \hfill
$\Box$

\bigskip

\noindent {\bf Proof of Proposition \ref{small-t}.}\  First we
claim that if $A$ is an arbitrary subset of $G$ of size at most
$\alpha n/8$, then $G-A$ contains a connected component of size at
least $\alpha n/4$. Indeed, if all components of $G-A$ have size
at most $\alpha n/4$, then by taking several of them together we
can find a subset $A'$ such that $\alpha n/4 \leq |A'| \leq \alpha
n/2$ and $A'$ has no neighbors in $G-A$, i.e., $N(A') \subseteq
A$.  On the other hand, by the remark in the beginning of the
section, we have that $|N(A')| \geq \alpha n/2$. This
contradiction proves our claim. Let
$$p=\frac{\alpha }{100}\sqrt{ \frac{ n \log t}{\log n}} \quad \mbox{and}\quad
q=\sqrt{ \frac{n \log n}{ \log t}},$$ and note that $pq=\alpha
n/100$. Hence , using the above  claim, we can greedily find $p$
disjoint sets $B_1, \ldots, B_p$, each of size $|B_i|=q$, such
that all induced subgraphs $G[B_i]$ are connected.

Let $B'=\cup_i B_i$, let $B''$ be an arbitrary subset of $G$ of size at most
$|B|'/10$ and let $B=B' \cup B''$. Then using the same argument as in the proof of Theorem
\ref{minor-t-expand} one can show that there exist a subset $X$ of $G-B$ of size at most $5|B|/t
\leq |B|/2$ such that the following
holds.
\begin{itemize}
\item
The graph $G'=G-X-B$ is a $(t/2, \alpha)$-expanding graph with at
most $\ell=2/\alpha$ connected components $G_1, \ldots, G_\ell$,
each of which therefore has diameter at most $7\alpha^{-1}\log
n/\log t$.
\item
There exists an index $1 \leq i \leq \ell$ such that at least
$p/(2\ell) \geq \alpha p/4$ sets $B_j$ have neighbors in $G_i$.
\end{itemize}
In particular this implies that there is a collection of $\alpha
p/4$ sets $B_j$, such that any pair of them can be connected by a
path $P$ of length at most $7\alpha^{-1}\log n/\log t$. Moreover
all vertices of $P$ except endpoints are contained in $G-B' \cup
B''$.

Now consider the following iterative procedure. In the beginning
of each iteration we will have sets $B'=\cup_i B_i$ and $B'',
|B''| \leq |B|'/10$, where $B''$ is the set of vertices of
disjoint paths that have been used at previous iterations to
connect sets $B_j$. We stop when we will have at least $\alpha
p/4$ sets $B_j$ which are pairwise connected. Then the contraction
of all these sets will give us a clique minor of size at least
$\Omega\left(\alpha^{2} \sqrt{\frac{n \log t}{\log n}}\right).$ By
the above discussion, at each iteration we indeed can construct a
path of length at most $7\alpha^{-1}\log n/\log t$ that does not
use vertices from $B' \cup B''$ and connects two previously not
connected sets $B_j$. Since the number of iterations is clearly at
most ${p \choose 2}$ we have that the size of the set $B''$
remains bounded by ${p \choose 2}7\alpha^{-1}\log n/\log t \leq
|B'|/10$ during all iterations. \hfill $\Box$

\section{Minors in pseudo-random graphs}\label{s-pseudo-random}
Here we prove Theorem \ref{jumbled-minor}. Throughout this section
we assume that $np$ is at least a sufficiently large constant and
$p$ is smaller than a sufficiently small constant.

\begin{lemma}
\label{jumbled-property1}
Let $G=(V,E)$ be a $(p,\beta)$-jumbled graph of order $n$ such that $ \beta=o(np)$.
Then $G$ contains an induced subgraph $G'$ of order $n'=(1-o(1))n$ such that
the degree of every vertex in $G'$ is $(1+o(1))n'p$ and every subset $X$ of $G'$
satisfies
$$e(X,V(G')-X)\geq (1-o(1))p|X|(n'-|X|).$$
\end{lemma}

\noindent {\bf Proof}\, Set $\epsilon=(4 \beta/(np))^{1/3}$ and
consider two disjoint subsets $S$ and $T$ both of size at least
$\epsilon n$. Then $e(S,T)=e(S\cup T)-e(S)-e(T)$ and therefore
\begin{eqnarray}
\label{cut}
e(S,T) &\geq& p{{|S|+|T|} \choose 2}-\beta(|S|+|T|)-p{|S| \choose 2}-\beta|S|-p{|T| \choose 2}-\beta |T|
\nonumber\\
&=& p|S||T|-2\beta(|S|+|T|) \geq p|S||T|-2\beta n \nonumber\\
&=& p|S||T|-\epsilon^3n^2p/2 \geq (1-\epsilon/2)p|S||T|.
\end{eqnarray}
Similarly one can show that
$e(S,T) \leq (1+\epsilon/2)p|S||T|$ for every two subsets $S, T$ as above.

Let $U$ be the set of vertices of $G$ with degree at least $(1+\epsilon)np$. If $U$ has
size at least $\epsilon n$ then we have that
\begin{eqnarray*}
e(U,V-U) &=& \sum_{v\in U} d(v)-2e(U)\geq (1+\epsilon)np|U|-2p{|U| \choose 2}-2\beta |U|\\
&\geq& (1+\epsilon)np|U|-p|U|^2-\epsilon^3np|U|\\
&>&(1+\epsilon/2)p|U|(n-|U|).
\end{eqnarray*}
This contradiction implies that there are less than $\epsilon n$
vertices in $G$ with degree at least $(1+\epsilon)np$. Let
$V_0=V-U$, $n_0=|V_0|>(1-\epsilon)n$,  and let $G_0$ be the
subgraph induced by $V_0$.

Consider the following process. If at step $i$ the graph $G_{i-1}$
contains a subset $X_i$ such that $|X_i|=x_i\leq \epsilon n$ and
$e(X_i,V(G_{i-1})-X_i)< (1-4\epsilon)px_i(n_{i-1}-x_i)$ delete
$X_i$ from the graph, update $G_i=G_{i-1}-X_i$, $n_i=|V(G_i)|$,
and continue. Consider the first time when we deleted at least
$\epsilon n$ vertices and let $Y=\cup_i X_i$. Then $\epsilon n
\leq |Y| \leq 2\epsilon n <3\epsilon n_0$ and
\begin{eqnarray*}
e(Y,V(G_0)-Y) &\leq& \sum_i e(X_i,V(G_{i-1})-X_i)<
(1-4\epsilon)p \sum_i x_i(n_{i-1}-x_i)\\
 &\leq& (1-4\epsilon)pn_0 \sum_i x_i = (1-4\epsilon)pn_0|Y|\\
&\leq&  \frac{1-4\epsilon}{1-3\epsilon} p|Y|(n_0-|Y|)
\leq (1-\epsilon/2)p|Y|(n_0-|Y|).
\end{eqnarray*}
This contradicts (\ref{cut}). Therefore there is a subset $Y$ of
$G_0$ of size at most $2\epsilon n$ such that every subset $X$ of
graph $G'=G_0-Y$ of size at most $\epsilon n$ satisfies
$e(X,V(G')-X) \geq (1-4\epsilon)p|X|(n'-|X|)$, where $n'=|V(G')|$.
In particular, taking $X$ to be a single vertex we have that the
minimum degree in $G'$ is at least $(1-4\epsilon)p(n'-1)$. By
(\ref{cut}) we also have that every subset $X$ with $\epsilon n
\leq |X| \leq n'/2$ satisfies that $e(X,V(G')-X) \geq
(1-4\epsilon)p|X|(n'-|X|)$. This inequality is satisfied by sets
of size larger than $n'/2$ by symmetry. Since $n'\geq
(1-3\epsilon)n$, by the above discussion, the maximum degree of
$G'$ is at most $(1+\epsilon)np \leq (1+5\epsilon)n'p$. Finally,
note that $\epsilon$ tends to zero as $np$ tends to infinity.
Therefore $G'$ satisfies the assertion of the lemma. \hfill $\Box$

\vspace{0.2cm}

A {\em lazy random walk} on a graph $G$ is a Markov chain whose
matrix of transition probabilities $P=(p_{i,j})$ is defined by
\begin{displaymath}
p_{i,j}= \left\{
\begin{array}{lll}
\frac{1}{2d(i)} & \textrm{if}~(i,j)\in E(G) \\
1/2& \textrm{if}~i=j\\
0 & \textrm{otherwise,}
\end{array} \right.
\end{displaymath}
i.e., if at some step we are at vertex $i$ than with probability
$1/2$ we stay at $i$ and with probability $\frac{1}{2d(i)}$ we
move to a random neighbor of $i$. This Markov chain has the
stationary distribution $\pi$ defined by
$\pi(i)=\frac{d(i)}{2e(G)}$. Let $\lambda_1\geq \lambda_2 \geq
\ldots \geq \lambda_n$ be the eigenvalues of $P$. Then the largest
eigenvalue $\lambda_1=1$ and since $P$ is positive semidefinite,
all other eigenvalues $\lambda_i, i \geq 2$, are non negative. For
more information about random walks on graphs we refer the
interested reader to the excellent survey of Lov\'asz \cite{L}.

\begin{lemma}
\label{jumbled-property2} Let $G=(V,E)$ be a graph of order $n$
such that every vertex in $G$ has degree $(1-o(1))np$ and  every
subset $X$ satisfies $e(X,V-X)\geq (1-o(1))p|X|(n-|X|)$. Then for
every subset $U$ of size $u$ the probability that a lazy random
walk on $G$ which starts from stationary distribution $\pi$ and
makes $\ell$ steps does not visit $U$ is at most
$e^{-0.03u\ell/n}$.
\end{lemma}

\noindent
{\bf Proof.}\, By the degree assumption we have that
$2|E|=\sum_id(i)=(1+o(1))n^2p$ and therefore the stationary distribution $\pi$
satisfies $\pi(i)=d(i)/(2|E|)=(1+o(1))/n$. Thus for every subset $S$ the measure of $S$ with respect to
$\pi$ equals $\pi(S)=\sum_{i \in S}\pi(i)=(1+o(1))|S|/n$.
Let
$$\Phi=\min_{\pi(S)\leq 1/2}\frac{\sum_{i \in S, j \in V-S}\,\, \pi(i)p_{i,j}}{\pi(S)\pi(V-S)},$$
be the {\em conductance} of $G$. By properties of $G$ we have that
\begin{eqnarray*}
\Phi&=&\min_{|S|\leq n/2+o(n)}(1+o(1))\frac{1}{n} \,\frac{1}{2np}\,\frac{e(S,V-S)}{(|S|/n)(1-|S|/n)}\\
&=&
\min_{|S|\leq n/2+o(n)}(1+o(1))\frac{1}{2n^2p}\frac{p|S|(n-|S|)}{(|S|/n)(1-|S|/n)}=1/2+o(1).
\end{eqnarray*}
Let $\lambda_2$ be the second largest eigenvalue of the transition
probabilities matrix $P$. Since all eigenvalues of $P$ are
non-negative we have that the {\em spectral gap} of this Markov
chain is $\delta=\max_{i\geq 2} 1-|\lambda_i|= 1-\lambda_2$. Then
by the result of Jerrum and Sinclair \cite{JS} (see also
\cite{L}), which provides a connection between the spectral gap
and the conductance of the graph, we have that
$\delta=1-\lambda_2\geq \Phi^2/8>0.03$. To finish the proof we can
now use well known estimates on the probability that a Markov
chain stays inside certain sets (see, e.g., \cite{AKS}, \cite
{AFWZ}, \cite{AS}, \cite{MORSS}). In particular, the assertion of
Theorem 5.4 in \cite{MORSS} implies that the probability that a
lazy random walk on $G$ which starts from stationary distribution
$\pi$ and makes $\ell$ steps does not visit a subset $U, |U|=u$ is
bounded from above by
$$\hspace{2cm}
\leq \big(1-\pi(U)\big)\big(1-\delta\pi(U)\big)^{\ell}
\leq \left(1-(1+o(1))\frac{\delta|U|}{n}\right)^{\ell} \leq e^{-0.03u\ell/n}. \hspace{2cm}\Box$$

\begin{lemma}
\label{jumbled-property3} Let $c>0$ be arbitrary constant. Let
$G=(V,E)$ be a $(p,\beta)$-jumbled graph of order $n$ such that $
\beta=o(np)$. Then $G$ contains a connected subset $B$ of size $c
p^{-1/2}$ such that it has at least $3c n\sqrt{p}/5$ neighbors in
$G$.
\end{lemma}

\noindent {\bf Proof}\, By Lemma \ref{jumbled-property1}, we can
assume that the minimum degree of $G$ is at least $(1+o(1))np$. We
construct $B$ using the following greedy procedure. Suppose we
have already constructed a connected set $B$ of size $k <c
p^{-1/2}$ which has at least $3knp/5$ neighbors in $G$. Let $X$ be
a subset of $3knp/5$ of these neighbors. Then the number of edges
inside $X\cup B$ is at most $p|X\cup B|^2/2+\beta|X\cup B
|<(np)|X|/6$. Therefore $X$ contains a vertex $v$ with at most
$np/3$ neighbors inside $X\cup B$. By the minimum degree
assumption $v$ more than $3np/5$ neighbors outside $X\cup B$.
Since by the definition of $X$, $v$ has also a neighbor in $B$,
the set $B\cup\{v\}$ is connected. This set has size $k+1$ and at
least $|X|+3np/5=3(k+1)np/5$ neighbors in $G$. Repeating this
process $c p^{-1/2}$ times we obtain a connected set $B$ that
satisfies the assertion of the lemma. \hfill $\Box$

\vspace{0.2cm} \noindent {\bf Proof of Theorem
\ref{jumbled-minor}.}\, Note that by definition any induced
subgraph of $G$ on at least $n/2$ vertices is still
$(p,2\beta)$-jumbled. Therefore by starting from $G$ and
repeatedly applying Lemma \ref{jumbled-property3} to the remaining
subgraph $G-\cup_{i<j}B_i$ we can construct $s=10^{-3}n\sqrt{p}$
disjoint connected sets $B_1, \ldots, B_s$ such that each $B_i$
has size $50p^{-1/2}$ and has at least $25n\sqrt{p}$ neighbors in
$G$. Let $D_1, \ldots, D_s$ be sets of size $25n\sqrt{p}$ such
that every vertex in $D_i$ has a neighbor in $B_i$. Consider the
following iterative procedure that we repeat $s$ times. In the
beginning of iteration $k+1$ we have connected sets $C_1, \ldots,
C_k$ each of size $50p^{-1/2}$, such that all $C_i$ and $B_j$ are
disjoint. We construct a new connected set $C_{k+1}$ of size
$50p^{-1/2}$ such that $C_{k+1}$ is disjoint from all previous
sets and there are at least $s/3$ indices $ 1 \leq j \leq s$ such
that there is an edge from $C_{k+1}$ to $B_j$. In the end of
algorithm if we contract all the sets $C_i, B_j$ we will get a
graph with average degree $\Omega(s)=\Omega(n\sqrt{p})$.

Let $U=(\cup_{i=1}^s B_i) \cup (\cup_{j \leq k}C_j)$ and note that
$|U| \leq n/10$. Then the induced subgraph $G[V\setminus U]$ is
$(p,2\beta)$-jumbled and therefore by Lemma
\ref{jumbled-property1} there is an induced subgraph $G'$ of $G-U$
on $n'\geq (1-o(1))(n-|U|) \geq 8n/9$ vertices such that the
degree of every vertex in $G'$ is $(1+o(1))n'p$ and every subset
$X$ of $G'$ satisfies $e(X,V(G')-X)\geq (1-o(1))p|X|(n'-|X|).$ Let
$V'$ be the vertex set of $G'$, $U'=V-V'$, and note that $|U'|
\leq n/9$. Next, we claim that $\sum_i|D_i-U'|\geq 10 n\sqrt{p}
\cdot s$. Note that from every vertex of $D_i \cap U'$ there is an
edge to one of the vertices in $B_i$. Since $B_i$ are disjoint,
each edge inside $U'$ is counted at most twice in the summation
$\sum_i|D_i-U'|$,  therefore $\sum_i |D_i \cap U'| \leq 2e(U')$.
This implied that
\begin{eqnarray*}
\sum_i |D_i-U'| &=& \sum_i \big(|D_i|-|D_i \cap U'|\big)\geq 25n\sqrt{p}s-e(U')\\
&\geq& 25n\sqrt{p}s-p|U'|^2/2-\beta|U'| \\
&\geq& 25n\sqrt{p}s -n^2p/80 -o(n^2p)\\
&\geq& 10 n\sqrt{p} s.
\end{eqnarray*}
Since  $|D_i-U'| \leq 25n\sqrt{p}$, we have that there are at
least $2s/5$ sets $D_i$ such that $D'_i=D_i-U'=D_i \cap V'$ has
size at least $10n\sqrt{p}$. Let $I$ be the set of indices $i$
such that $|D'_i| \geq 10n\sqrt{p}$.

Consider a lazy random walk on $G'$ which starts from the
stationary distribution and makes $\ell=50p^{-1/2}$ steps. By
Lemma \ref{jumbled-property2} the probability that this walk does
not intersects a given $D'_i, i\in I$, is at most
$e^{-0.03|D_i'|\ell/n'} \leq 0.01$. Therefore by Markov's
inequality with positive probability this walk intersects at least
$0.9|I|\geq s/3$ sets $D'_i$. Choose one such walk and denote its
vertex set by $C_{k+1}$. This gives a connected subset of size (at
most)  $50p^{-1/2}$, which by definition is disjoint from all
previous sets $B_i, C_j$ and has neighbors in at least $s/3$ sets
$B_i$. \hfill $\Box$

\section{Minors in $H$-free graphs}\label{s-H-free}
In this section we prove Theorems \ref{K_{s,s'}-minor} and
\ref{H-minor}.

We start with proving Theorem \ref{K_{s,s'}-minor}. We assume that
$s,s'$ are fixed integers satisfying $2\le s\le s'$.

\begin{lemma}
\label{small-separation} Let $G$ be a graph of order $n$ with
average degree $d \leq r$. Let $X, Y, Z$ be a partition of the
vertex set of $G$ into three disjoint sets such that $|Y| \leq
\frac{|X|}{2a}$ and $e(X,Z)\leq \frac{r}{4a}|X|$ for some $a>0$.
Then $G\setminus X$ still has the average degree at least $d$, or
the average degree of the subgraph induced by the set $X\cup Y$ is
at least $d-\frac{r}{a}$.
\end{lemma}

\noindent {\bf Proof.}\, Let $|X|=\alpha n$ and suppose that the
average degree of $G \setminus X$ is at most $d$, i.e.,
$e(G\setminus X) \leq (1-\alpha)dn/2$. Let $G'$ be the subgraph of
$G$ induced by the set $X\cup Y$. Then $|V(G')|=|X\cup Y|\leq
(1+1/(2a))\alpha n$ and
\begin{eqnarray*}
e(G')&\geq& e(G)-e(G \setminus X)-e(X,Z) \\
&\geq&
dn/2-(1-\alpha)dn/2-\frac{r}{4a}\alpha n\\
&=&\Big(d-\frac{r}{2a}\Big)\frac{\alpha n}{2}.
\end{eqnarray*}
Since $d \leq r$, the average degree of $G'$ is:
$$\hspace{3.5cm}
\frac{2e(G')}{|V(G')|}\geq \frac{(d-r/(2a))\alpha n}{(1+1/(2a))\alpha n}=
\frac{2a}{2a+1}d-\frac{r}{2a+1}
\geq d-\frac{r}{a}. \hspace{3.5cm} \Box$$

\begin{lemma}
\label{K_{s,s'}-expand}
Let $G$ be $K_{s,s'}$-free graph, $s'\geq s$ and let $X \subseteq V(G)$ such
that $e(X,V-X) \geq d|X|$ for some $d>0$. Then
\begin{displaymath}
|N(X)| \geq \left\{
\begin{array}{ll}
\frac{d|X|}{s'} & \textrm{if}~|X| \leq d^{1/(s-1)} \\
$~$&  $~$\\
\frac{d^{s/(s-1)}}{s'} & \textrm{otherwise}
\end{array} \right.
\end{displaymath}
\end{lemma}

\noindent {\bf Proof.}\, First note that we need only to consider
the case when $|X| \leq d^{1/(s-1)}$. Indeed if $|X|\geq
d^{1/(s-1)}$ then by the averaging argument there exists
$X'\subseteq X$ of size $|X'|=d^{1/(s-1)}$ such that $e(X',V-X)
\geq d|X'|$.

Let $|X|\le d^{1/(s-1)}$. Assume by the way of contradiction that
$|N(X)|< d|X|/s'$. Let $Y$ be a subset of $d|X|/s'$ vertices of
$V\setminus X$ containing $N(X)$. Then there are at least $d|X|$
edges between $X$ and $Y$ in $G$. Let us count the number of pairs
$(y,S)$, where $y\in Y$, $S\subseteq X\cap N(y)$, $|S|=s$. Denote
this quantity by $A$. Then
$$
A= \sum_{y\in Y}{{d(y,X)}\choose s} \ge
|Y|{{\frac{\sum_{y\in Y} d(y,X)}{|Y|}}\choose s} \ge
|Y|
{{\frac{d|X|}{|Y|}}\choose s}=\frac{d|X|}{s'}{{s'}\choose s}\ .
$$
On the other hand, each $S$ appears in at most $s'-1$ pairs
$(y,S)$ as otherwise we get a copy of $K_{s,s'}$ with $s$ vertices
in $S$ and $s'$  vertices in $X$. Therefore,
$$
A\le (s'-1){{|X|}\choose s}\ .
$$
Comparing the above two estimates for $A$ we get:
$$
\frac{d|X|}{s'}{{s'}\choose s}\le A\le (s'-1){{|X|}\choose s}<
(s'-1)\frac{|X|^s}{s!}\,,
$$
implying:
$$
\frac{s!}{s'(s'-1)}{{s'}\choose s} < \frac{|X|^{s-1}}{d}\ .
$$
As $s'\ge 2$, the LHS of the inequality above is easily seen to be
at least 1, while by the assumption $|X|\le d^{1/(s-1)}$, the RHS
is at most 1 -- a contradiction. \hfill $\Box$

\begin{lemma}
\label{l45}
Let $c>0$ be a constant and let
$G$ be a $K_{s,s'}$-free graph on $cr^{\frac{s}{s-1}}$
vertices with average degree $r$.
Then $G$ contains a minor with average degree at least
$\Omega\big(r^{1+\frac{1}{2(s-1)}}\big)$.
\end{lemma}

\noindent {\bf Proof.}\, Since the average degree of $G$ is at
least $r$, it contains a subgraph $G'$ with minimum degree at
least  $r/2$. Let $X$ be a subset of $G'$ of size at most
$r^{\frac{1}{s-1}}/4$. Since the minimum degree is at least $r/2$,
every vertex of $X$ has at least $r/4$ neighbors outside $X$,
i.e., $e_{G'}(X,V(G')-X) \geq \frac{r}{4}|X|$. Therefore by Lemma
\ref{K_{s,s'}-expand} we have that $|N_{G'}(X)| \geq
\frac{r}{4s'}|X|$. This implies that $G'$ is a $(t,
\alpha)$-expanding graph of order $n= cr^{\frac{s}{s-1}}$, where
$t=r/(4s')$ and $ \alpha=\frac{1}{16s'c}$. Thus, by Theorem
\ref{minor-t-expand}, it contains a minor with average degree at
least
$$\hspace{5.2cm} \Omega\left(\alpha^3 \frac{\sqrt{nt \log t}}{\sqrt{\log n}}\right)=
\Omega\left(r^{1+\frac{1}{2(s-1)}}\right).
\hspace{5.2cm} \Box
$$

\begin{lemma}
\label{diameter-K_{s,s'}} Let $2 \leq s \leq s' \leq a$ and let
$G$ be a $K_{s,s'}$-free graph of order $n \leq
e^{2a}r^{\frac{s}{s-1}}$ such that for every two disjoint subsets
$X, |X|\leq n/2$, and $Y,|Y|\leq \frac{1}{3a^2}|X|$, we have that
$e\big(X,V(G)-(X\cup Y)\big) \geq \frac{r}{4a^2}|X|$. Then the
diameter of $G$ is at most $33a^3$.
\end{lemma}

\noindent {\bf Proof.}\, By the above condition, $G$ has minimum
degree  at least $\frac{r}{4a^2}$. If $\frac{r}{4a^2}>n/2$ we are
done, since the diameter of $G$ is at most two. Let $v$ be an
arbitrary vertex of $G$ and let $X \subset N(v)$ be a subset of
$\frac{r}{4a^2}$ neighbors of $v$. Our assumptions on $G$ imply
that $e_G(X,V-X) \geq \frac{r}{4a^2}|X|$. Since $G$ is
$K_{s,s'}$-free, $s \geq 2$ and $s' \leq a$, by Lemma
\ref{K_{s,s'}-expand} (with $d=\frac{r}{4a^2})$, we have that
$$|N(X)| \geq \min\left\{\frac{r}{4a^2s'}|X|,
\frac{1}{s'}\Big(\frac{r}{4a^2}\Big)^{\frac{s}{s-1}}\right\}\geq
\frac{r^{\frac{s}{s-1}}}{16a^5}.
$$
Therefore there are at least $\frac{1}{16a^5}r^{\frac{s}{s-1}}$
vertices within distance at most two from any vertex of $G$. We
also have  that every subset $U$ of $G$ of size at most $n/2$
satisfies $|U \cup N(U)| \geq \big(1+\frac{1}{3a^2}\big)|U|$.
Since $8a^5e^{2a} < \big(1+\frac{1}{3a^2}\big)^{16a^3}$ for $a\ge
2$ , we conclude that there are more than
$$\frac{r^{\frac{s}{s-1}}}{16a^5}\left(1+\frac{1}{3a^2}\right)^{16a^3} >
\frac{1}{2}e^{2a} r^{\frac{s}{s-1}}\ge n/2$$ vertices within
distance at most $2+16a^3$ from any given vertex of $G$. This
implies that the diameter of $G$ is at most $2(2+16a^3) \leq
33a^3$. \hfill $\Box$

\begin{lemma}
\label{large-minor} Let $2 \leq s \leq s' \leq a \leq 2\log r$ and
let $G$ be a $K_{s,s'}$-free graph of order $n$ such that
$a^{14}r^{\frac{s}{s-1}}  \leq n \leq e^{2a}r^{\frac{s}{s-1}}$ and
for every two disjoint subsets $X, |X|\leq 0.7n$ and $Y,|Y|\leq
\frac{1}{2a^2}|X|$ we have that $e\big(X,V(G)-(X\cup Y)\big) \geq
\frac{r}{4a^2}|X|$. Then $G$ contains a minor with average degree
at least $cr^{1+\frac{1}{2(s-1)}}$, where $c>0$ is a constant
independent of $r$ and $a$.
\end{lemma}

\noindent
{\bf Proof.} \,
Let
$$p=\frac{1}{10^3}a^2r^{1+\frac{1}{2(s-1)}} \quad \mbox{and}\quad
q=33 \frac{n}{a^4r^{1+\frac{1}{2(s-1)}}},$$ and consider the
following iterative procedure which we will repeat $p$ times. In
the beginning of iteration $k+1$ we will have $k$ disjoint sets
$B_1, \ldots, B_k$ each of size $|B_i|=q$ such that all induced
subgraph $G[B_i]$ are connected. We will construct a new subset
$B_{k+1}$, also of size $q$, such that the induced subgraph
$G[B_{k+1}]$ is connected and there are at least $k/(8a^2)$
indices $1 \leq i\leq k$ such that there is an edge from $B_i$ to
$B_{k+1}$. In the end of this algorithm, if we contract all
subsets $B_i$ we will get a graph with average degree
$\Omega(\frac{p}{8a^2})\geq c r^{1+\frac{1}{2(s-1)}}$.

Let $B=\cup_{i=1}^k B_i$ and note that $|B|\leq pq \leq
\frac{n}{30a^2}$. Denote $C=V(G)-B$ and let $G'$ be the subgraph
of $G$ induced by $C$. Let $X_1$ and $Y_1$ be two disjoint subsets
of $C$ such that $n/5 \leq |X_1|\leq 0.7n$, $|Y_1| \leq
\frac{1}{3a^2}|X_1|$ and $e\big(X_1,C-(X_1\cup Y_1)\big) <
\frac{r}{4a^2}|X_1|$. Set  $Y'=Y_1 \cup B$. Then we have
$$|Y'|\leq |Y_1|+|B| \leq \frac{1}{3a^2}|X_1|+\frac{n}{30a^2} \leq \frac{1}{2a^2}|X_1|$$
and  $e\big(X_1,V(G)-(X_1\cup Y')\big) < \frac{r}{4a^2}|X|$ which
contradicts our assumption about $G$. Therefore there exist two
disjoint (or empty) subsets $X_1,Y_1 \subset C$ such that $|X_1|
\leq n/5$, $|Y_1| \le \frac{1}{3a^2}|X_1|$, $e\big(X_1,C-(X_1\cup
Y_1)\big) \le \frac{r}{4a^2}|X_1|$  and the remaining set
$D=C-X_1$ spans a graph $G''$ in which for every two disjoint
subsets $X, |X|\leq n/2$, and $Y,|Y|\leq \frac{1}{3a^2}|X|$, we
have that $e\big(X,V(G'')-(X\cup Y)\big) \geq \frac{r}{4a^2}|X|$.
(Such a pair $(X_1,Y_1)$ can be obtained by repeatedly deleting
sets $(X,Y)$ such that $0<|X|\le n/5$, $|Y|\le |X|/(3a^2)$ and in
the obtained graph $G''$, $e(X,V(G'')-(X\cup Y))\le
\frac{r}{4a^2}$, for as long as the union of the deleted $X$'s
does not go over $n/5$.) Note that by Lemma
\ref{diameter-K_{s,s'}}, $G''$ has diameter at most $33a^3$.

Consider all sets $B_j$ that satisfy $e(B_j,D) \geq \frac{r}{4a^2}
|B_j|$. Without loss of generality, we can assume that the first
$m$ sets $B_1, \ldots, B_m$ have this property. We claim that $m$
is at least $\frac{k}{4a^2}$. If this is not the case then denote
$Y_2=\cup_{j=1}^m B_j$, and $X_2=\cup_{j=m+1}^k B_j$. By
definition $|Y_2| \leq \frac{m}{k-m}|X_2|\leq \frac{1}{3a^2}|X_2|$
and
$$e(X_2,D)=\sum_{j=m+1}^ke(B_j,D) < \sum_{j=m+1}^k \frac{r}{4a^2} |B_j|=
\frac{r}{4a^2}|X_2|.$$
Define $X=X_1 \cup X_2$ and $Y=Y_1 \cup Y_2$. Then $|X| \leq n/5+|B|\leq n/4$,
$$|Y| \leq |Y_1|+|Y_2| \leq \frac{1}{3a^2}|X_1|+\frac{1}{3a^2}|X_2| \leq \frac{1}{2a^2}|X|,$$
and also
\begin{eqnarray*}
e\big(X,V(G)-(X\cup Y)\big) &\leq& e(X_1,D-Y_1)+e(X_2,D)\\
&<&\frac{r}{4a^2}\big(|X_1|+|X_2|\big)
=\frac{r}{4a^2}|X|.
\end{eqnarray*}
This contradicts the properties of $G$. Therefore we have that the
first $m=\frac{k}{3a^2}$ sets $B_1, \ldots, B_m$ satisfy that
$e(B_j,D) \geq \frac{r}{4a^2} |B_j|$.

Denote by $U_j, 1 \leq j \leq m$, the set of neighbors of $B_j$ in
$D$. Since $n \geq a^{14} r^{\frac{s}{s-1}}$, $s' \leq a$,
$|B_j|=q$ and $a=r^{o(1)}$, By Lemma \ref{K_{s,s'}-expand} (with
$d=\frac{r}{4a^2})$, we have that
$$|U_j| \geq \min\left\{\frac{r}{4a^2s'}|B_j|,
\frac{1}{s'}\Big(\frac{r}{4a^2}\Big)^{\frac{s}{s-1}}\right\} \geq
a^7r^{1+\frac{1}{2(s-1)}}.$$
 Pick uniformly at random with
repetition $n/(a^7r^{1+\frac{1}{2(s-1)}})$ vertices of $G''$
 and denote this set by $W$. For every index $1 \leq i
\leq m$ the probability that $W$ does not intersect $U_i$ is at
most $\left(1-\frac{|U_i|}{|G''|}\right)^{|W|}\leq 1/e$. Therefore
the expected number of sets $U_i$ which have non-empty
intersection with $W$ is at least $(1-1/e)m>m/2$. Hence there is a
choice of $W$ that intersects at least $m/2\geq k/(8a^2)$ sets
$U_i$. Fix an arbitrary vertex  $w_0 \in W$ and consider a
collection of shortest paths in $G_1$ from $w_0$ to the remaining
vertices in $W$. Since the diameter of $G''$ is at most $33a^3$
and $33a^3|W| \leq q$, by taking union of these paths and adding
extra vertices if necessary we can construct a connected subset of
size $q$ containing $W$. Denote this set by $B_{k+1}$ and note
that it is connected by an edge to at least $k/(8a^2)$ sets $U_i,
i \leq k$. This completes the proof of the lemma. \hfill $\Box$

\begin{lemma}
\label{smallgraphs} Let $G$ be a $K_{s,s'}$-free graph of average
degree $r$ and at most $r^{4+\frac{s}{s-1}}$ vertices. Then $G$
contains a minor with average degree at least
$$\Omega\left(r^{1+\frac{1}{2(s-1)}}\right).$$
\end{lemma}

\noindent {\bf Proof.}\, Let $\{a_i, i \geq 0\}$ be an increasing
sequence defined by $a_0=20s'$ and $a_{i+1}=e^{a_i/7}$. Note that
$a^{14}_{i+1}=e^{2a_i}$ and let $\ell$ be the first index such
that $e^{2a_{\ell}} > r^{4+\frac{s}{s-1}}$. Then there is  some $0
\leq i \leq \ell$ so that the order $n$ of our graph $G$ satisfies
$$a_i^{14}r^{\frac{s}{s-1}}  \leq n < e^{2a_i}r^{\frac{s}{s-1}}.$$
If $G$ has the property that for every two disjoint subsets $X,
|X|\leq 0.7n$, and $Y,|Y|\leq \frac{1}{2a_i^2}|X|$, we have that
$$e\big(X,V(G)-(X\cup Y)\big) \geq \frac{r}{4a_i^2}|X|,$$ then by
Lemma \ref{large-minor} it contains a minor with average degree
$\Omega\left(r^{1+\frac{1}{2(s-1)}}\right)$ and we are done.
Otherwise, there are two sets $X, Y$ as above for which
$e\big(X,V(G)-(X\cup Y)\big) < \frac{r}{4a_i^2}|X|.$ Then, by
Lemma \ref{small-separation}, we either have that the average
degree of graph $G-X$ is at least $r$, or the average degree of
the subgraph induced by $X \cup Y$ is at least
$r-\frac{r}{a_i^2}$. In the first case let $G_1=G-X$ and in the
second let $G_1=G[X \cup Y]$. Note that the number of vertices
$n_1$ of new graph is strictly smaller than that of $G$. Moreover
if the average degree of $G_1$ is smaller than that of $G$ we know
that $n_1=|X\cup Y|\leq 3n/4$. Continue this process until we
either find a minor  with average degree at least
$\Omega\left(r^{1+\frac{1}{2(s-1)}}\right)$, or arrive to a graph
$G'$ with $n'$ vertices such that $n' \leq
a_0^{14}r^{\frac{s}{s-1}}$.

In the first case we are obviously done. In the second case we
claim that the average degree of $G'$ is still at least $r/2$.
Note that if at some stage the order of our graph $G_j$ satisfied
$$ a_i^{14}r^{\frac{s}{s-1}}  < |V(G_j)| \leq e^{2a_i}r^{\frac{s}{s-1}},$$
then the average degree of the new graph $G_{j+1}$ could decrease
only by at most $r/a_i^2$. In this case the order of $G_{j+1}$
drops as well so that $|G_{j+1}| \leq 3|G_j|/4$. Since
$(3/4)^4<e^{-1}$, we have that this can happen only at most $8a_i$
times, before the size of the remaining graph will become smaller
than $a_i^{14}r^{\frac{s}{s-1}}=e^{2a_{i-1}}r^{\frac{s}{s-1}}$.
Since $a_{i+1}=e^{a_i/7}\geq 2a_i$,  we have that during all
iterations the average degree of the resulting graph can decrease
by at most
$$ \sum_i8a_i \cdot \frac{r}{a_i^2}=r \sum_i \frac{8}{a_i} \leq \frac{16}{a_0} r <r/2.$$
Hence the final graph $G'$ has average degree at least $r/2$ and at most
$O\big(r^{\frac{s}{s-1}}\big)$ vertices. Therefore, by Lemma \ref{l45}, it contains a
minor with average degree $\Omega\left(r^{1+\frac{1}{2(s-1)}}\right)$.
\hfill $\Box$

\medskip

 \noindent {\bf Proof of Theorem \ref{K_{s,s'}-minor}.}\,
Let $G$ be a $K_{s,s'}$-free graph with average degree $r$ and let
$n$ be the number of vertices of $G$. By Lemma \ref{smallgraphs},
we can assume that $n > r^5$. Suppose that $G$ contains a subset
$X, |X| \leq 0.7n$, such that $|N(X)| \leq \frac{|X|}{2\log^2 n}$.
If the average degree of $G-X$ is at least $r$, set $G_1=G-X$ and
let $n_1$ be the number of vertices in $G_1$. Otherwise, let $G_1$
be the subgraph induced by the set $X\cup N(X)$. In the second
case, by Lemma \ref{small-separation}, the average degree of $G_1$
is at least $r-\frac{r}{\log^2 n}$. Note that in both cases we
obtain a smaller graph. Moreover if the average degree of $G_1$ is
smaller than that of $G$ we know that $n_1=|X\cup N(X)|\leq 3n/4$.
Continue this process until we obtain a subgraph $G'$ of $G$ on
$n'$ vertices such that one of the following holds. Either $n'
\leq r^5$ or every subset $X$ of $G'$ of size $|X| \leq 0.7n'$ has
$|N(X)|  \geq \frac{|X|}{2\log^2 n'}$. Note that in the second
case the graph $G'$ does not have a separator of size
$\frac{n'}{2\log^2 n'}$. Since $n'>r^5$, by a result of Plotkin,
Rao and Smith \cite{PRS}, $G'$ has a clique minor of size
$$\Omega\left(\frac{n'/\log^2 n'}{\sqrt{n' \log n'}}\right)
=\Omega\left(\frac{\sqrt{n'}}{\log^{5/2} n'}\right) \geq
r^{5/2-o(1)} \gg r^{1+\frac{1}{2(s-1)}}.$$ In the first case, when
$n' \leq r^5$ we claim that the average degree of $G'$ is still at
least $r/2$. Indeed, let $x_0=n, x_1, \ldots, x_\ell \geq r^5$ be
the sequence of orders of graphs that  we had during the process
when the average degree decreased. Then we know that $x_{i+1}\leq
3 x_i/4$ and the decrease in the average degree at the
corresponding step was at most $r/\log^2 x_i$. Let $y_i=\log
x_{\ell-i}$, then $y_0 \geq 5\log r$ and $y_{i+1}\geq
y_i+\log(4/3) \geq y_i+1/4$. Therefore
$$ \sum_i \frac{1}{y^2_i} \leq \sum_i \frac{1}{(y_0+i/4)^2} \leq 16 \sum_{i=0}^{\infty}
\frac{1}{(4y_0+i)^2} \leq \frac{4}{y_0-1}\ll 1/2,$$ and we
conclude that the average degree of $G'$ is at least
$r\big(1-\sum_i 1/\log^2 x_i\big) \geq r/2$. Therefore we can find
in $G'$ a minor with average degree
$\Omega\left(r^{1+\frac{1}{2(s-1)}}\right)$ using Lemma
\ref{smallgraphs}. This completes the proof of the theorem. \hfill
$\Box$

\bigskip

The proof of Theorem \ref{H-minor} is very similar to that of
Theorem \ref{K_{s,s'}-minor}. The only (relatively) substantial
difference in the proof of Theorem \ref{H-minor} compared to that
of \ref{K_{s,s'}-minor} lies in the proof of Lemma
\ref{K_{s,s'}-expand}. Instead, we have:

\begin{lemma}
\label{H-expand}
Let $H$ be a bipartite graph of order $h$ with parts $A$ and $B$
such that the degrees of all vertices
in $B$ do not exceed $s$.
Let $G$ be $H$-free graph and let $X \subseteq V(G)$ such that
$e(X,V-X) \geq (2dh)|X|$ for some $d>0$. Then
\begin{displaymath}
|N(X)| \geq \left\{
\begin{array}{ll}
d|X| & \textrm{if}~|X| \leq d^{1/(s-1)} \\
d^{s/(s-1)} & \textrm{otherwise}
\end{array} \right.
\end{displaymath}
\end{lemma}

\noindent {\bf Proof.}\, Similarly to the proof of Lemma
\ref{K_{s,s'}-expand} we need to consider only the case when $|X|
\leq d^{1/(s-1)}$. Then the result follows from a variant of the
dependent random choice argument utilized in particular in
\cite{AKS2}. If $|N(X)|\leq d|X|$ then pick a random vertex $v$ in
$N(X)$. Let the random variable $Y$ count the number of neighbors
of $v$ in $X$, and let the random variable $Z$ be the number of
$s$-tuples of vertices in $N(v)\cap X$ that have at most $s-1$
common neighbors. Then the expected value of $Y$ is at least
$\frac{e(X,V-X)}{|N(X)|} \geq \frac{2dh|X|}{|N(X)|}$, while the expected value of $Z$ is at most
${{|X|}\choose s}\frac{h-1}{|N(X)|}$. It thus follows that
\begin{eqnarray*}
 \mathbb{E}[Y-Z] &=& \mathbb{E}[Y]-\mathbb{E}[Z]\ge \frac{2dh|X|}{|N(X)|}-{{|X|}\choose
s}\frac{h-1}{|N(X)|} \\
&>& \frac{h}{|N(X)|}\left(2d|X|-{{|X|}\choose
s}\right)\\
&>& \frac{h|X|}{|N(X)|}\left(2d-\frac{|X|^{s-1}}{s!}\right) \\
&\ge&
\frac{dh|X|}{|N(X)|}
\ge h\ .
\end{eqnarray*}

Therefore there exists a vertex $v\in N(X)$ so that
$Y-Z\ge h$. Fix such a vertex, denote by $A_0$ its neighborhood in
$X$, and for each $s$-tuple $S$ in $A_0$ with less than $h$ common
neighbors, delete an arbitrary vertex from $S$. Denote the
obtained set by $A_1$. Then $|A_1|\ge h$, and every $s$-tuple in
$A_1$ has at least $h$ common neighbors. We then can embed a copy
of $H$ in $G$ by first embedding the side $A$ of $H$ one-to-one
into $A_1$, and then embedding the vertices of $B$, the other side
of $H$, vertex by vertex. As every $s$-tuple in $A_1$ has at least
$h$ common neighbors and the degree of every vertex in $B$ is at most $s$, we will be always able
to find a required
vertex.
 \hfill $\Box$

\medskip

Repeating the proof of Theorem \ref{K_{s,s'}-minor} and using the
above lemma instead of Lemma \ref{K_{s,s'}-expand}, we can prove
Theorem \ref{H-minor}.

\section{Minors in $C_{2k}$-free graphs}\label{s-C_2k-free}
Here we prove Theorem \ref{C_2k-minor}. In the rest of this
section we may and will assume that $k\geq 3$ is fixed ($k=2$
follows from Theorem \ref{K_{s,s'}-minor}) and $r$ is sufficiently
large compared to $k$.

\begin{lemma}
\label{C_2k-size}
Let $G$ be $C_{2k}$-free graph on $n$ vertices with average degree $d$.
Then $n \geq \big(\frac{d}{16k}\big)^k$.
\end{lemma}

\noindent {\bf Proof.}\, It was proved in \cite{V} that the number
of edges in a $C_{2k}$-free graph on $n$ vertices is at most  $8k
n^{1 + \frac{1}{k}}$. Therefore we have that $nd/2 \leq 8k n^{1 +
\frac{1}{k}}$, which implies that $n \geq
\big(\frac{d}{16k}\big)^k$. \hfill $\Box$

\begin{lemma}
\label{C_2k-expand} Let $k \geq 3$ and let $G$ be a $C_{2k}$-free
graph. If $X \subseteq V(G)$ satisfies that $e(X,V-X) \geq d|X|$
for some $d \geq k$, then
\begin{displaymath}
|N(X)| \geq \left\{
\begin{array}{lll}
\frac{d|X|}{4k^2} & \textrm{if}~|X| \leq d^{\frac{k-1}{2}} \\
$~$&  $~$\\
\frac{d^{1/2}|X|}{4k^2} & \textrm{if}~|X| \leq d^{\frac{k+1}{2}} \\
$~$&  $~$\\
3|X| & \textrm{if}~|X| \leq \Big(\frac{d}{6k}\Big)^k
\end{array} \right.
\end{displaymath}
\end{lemma}

\noindent {\bf Proof.}\, This estimates can be easily deduced from
a result of Naor and Verstra\"{e}te \cite{NV}, who proved that the
number of edges in a $C_{2k}$-free bipartite graph with parts $X$
and $Y$ is bounded by
$$e(X,Y) \leq (2k-3)\Big((|X||Y|)^{\frac{k+1}{2k}}+|X|+|Y|\Big).$$
Indeed, we will have a contradiction with this inequality if $e(X,N(X)) \geq d|X|$ and the size of $N(X)$ is
less than in the assertion of the lemma.
\hfill $\Box$

\begin{lemma}
\label{diameter-C_2k} Let $k \geq 3, \alpha \geq 1, \rho \geq 3$
and let $G$ be a $C_{2k}$-free graph of order $n \leq \rho r^k$
such that for every two disjoint subsets $X, |X|\leq n/2$, and
$Y,|Y|\leq \frac{1}{3\alpha}|X|$, we have that
$e\big(X,V(G)-(X\cup Y)\big) \geq \frac{r}{4\alpha}|X|$. Then
every subset $W \subset G$ of size at least $r^{k/2-1}\log r$ is
contained in a connected subgraph of $G$ on at most
$(40k^2\alpha^{3/2} \log \rho)|W|$ vertices.
\end{lemma}

\noindent {\bf Proof.}\, By the above condition and Lemma
\ref{C_2k-expand}, $G$ has minimum degree  at least
$\frac{r}{4\alpha}$ and every subset of $G$ of size at most
$\big(\frac{r}{24k \alpha}\big)^k$ expands at least three times.
Therefore for every vertex $v$ there are at least
$\big(\frac{r}{24k \alpha}\big)^k$ vertices which are within
distance at most $k \log r$ from $v$. We also have  that every
subset $U$ of $G$ of size at most $n/2$ satisfies $|U \cup N(U)|
\geq \big(1+\frac{1}{3\alpha}\big)|U|$. Since $\rho \big(24k
\alpha\big)^k < \big(1+\frac{1}{3\alpha}\big)^{4\alpha \log
\rho+8k^2\alpha^{3/2} }$, we conclude that there are more than
$$\frac{r^k}{(24k \alpha)^k}\left(1+\frac{1}{3\alpha}\right)^{4\alpha \log \rho+8k^2\alpha^{3/2}} >
\frac{1}{2} \rho r^{k}=n/2$$ vertices within distance at most $k
\log r+4\alpha \log \rho+8k^2\alpha^{3/2}$ from any given vertex
of $G$. This implies that the diameter of $G$ is at most $2k\log
r+8\alpha \log \rho +16k^2\alpha^{3/2}$.

Similarly, by Lemma \ref{C_2k-expand}, there are at least
$\big(\frac{r}{16k^2 \alpha}\big)^{\frac{k+1}{2}}$ vertices within
distance at most $\frac{k+1}{2}$ from every vertex of $G$, and
therefore the number of vertices within distance at most
$\frac{k+1}{2}+1 \leq k$ is at least
$$\frac{(r/(4\alpha))^{1/2}}{4k^2}\left(\frac{r}{16k^2\alpha}\right)^{\frac{k+1}{2}} \geq
\frac{r^{k/2+1}}{(16k^2\alpha)^{k}}.$$
Since $\rho \big(16k^2\alpha\big)^{k} < \big(1+\frac{1}{3\alpha}\big)^{4\alpha \log \rho
+8k^2\alpha^{3/2}}$,
we conclude that there are more than
$$\frac{r^{k/2+1}}{(16k^2\alpha)^{k}}\left(1+\frac{1}{3\alpha}\right)^{4\alpha \log \rho
+8k^2\alpha^{3/2}} > \rho r^{k/2+1}$$ vertices within distance at
most $k+4\alpha \log \rho +8k^2\alpha^{3/2} \leq 4\alpha \log \rho
+9k^2\alpha^{3/2}$ from any given vertex of $G$.

Let $W$ be a subset of $V(G)$ of size at least $r^{\frac{k-1}{2}}$
and consider the following iterative process that constructs a
connected subgraph $G'$ of $G$ containing $W$. At the beginning
the vertex set of $G'$ is $W$. At every step if there are two
connected components of $G'$ such that the distance between them
is at most $8\alpha \log \rho +18k^2\alpha^{3/2}$, connect them by
a shortest path and add the vertices of this path to $G'$. We
perform this step at most $|W|$ times until the distance between
every two remaining connected components of $G'$ is larger than
$8\alpha \log \rho +18k^2\alpha^{3/2}$. Then the balls of radius
$4\alpha \log \rho +9k^2\alpha^{3/2}$ around each component are
disjoint. By the above discussion, each such ball contains at
least $\rho r^{k/2+1}$ vertices so the number of components is at
most $\frac{n}{\rho r^{k/2+1}} \leq r^{k/2-1}$. Now fix one
component of $G'$ and connect it to every other component by a
path whose length is bounded by the diameter of $G$. This gives a
connected subgraph of $G$ that contains $W$ and has altogether at
most
$$ \big(8\alpha \log \rho +18k^2\alpha^{3/2})|W|+r^{k/2-1}
(2k\log r+8\alpha \log \rho +16k^2\alpha^{3/2}) \leq \big(40 k^2\alpha^{3/2} \log \rho)|W|$$
vertices. This completes the proof of the lemma.
\hfill $\Box$

\begin{lemma}
\label{large-minor-C_2k} Let $\alpha \geq 1,  3 \leq \rho \leq
r^2$, and let $G$ be a $C_{2k}$-free graph of order $n \leq \rho
r^k$ such that for every two disjoint subsets $X, |X|\leq 0.7n$,
and $Y,|Y|\leq \frac{1}{2\alpha}|X|$, we have that
$e\big(X,V(G)-(X\cup Y)\big) \geq \frac{r}{4\alpha}|X|$. Then $G$
contains a minor with average degree at least
$$c \frac{r \cdot n^{\frac{k-1}{2k}}}{\alpha^{\frac{15k+3}{4k}} \log^{\frac{k+1}{2k}} \rho},$$
where $c$ is a constant independent of $r, \rho$
and $\alpha$.
\end{lemma}

\noindent
{\bf Proof.} \,
Let
$$q=\frac{100k^3 \alpha}{r}\big(n \alpha^{3/2}\log \rho\big)^{\frac{k+1}{2k}}
\quad \mbox{and}\quad p=\frac{n}{30\alpha q},$$ and consider the
following iterative procedure which we will repeat $p$ times. In
the beginning of iteration $t+1$ we  have $t$ disjoint sets $B_1,
\ldots, B_t$, each of size $|B_i|=q$, such that all induced
subgraphs $G[B_i]$ are connected. We will construct a new subset
$B_{t+1}$, also of size $q$, such that the induced subgraph
$G[B_{t+1}]$ is connected, and there are at least $t/(8\alpha)$
indices $1 \leq i\leq t$ such that there is an edge from $B_i$ to
$B_{t+1}$. In the end of this algorithm if we contract all subsets
$B_i$ we get a graph with average degree
$$\Omega\Big(\frac{p}{8\alpha}\Big)\geq \Omega\left(
\frac{r \cdot n^{\frac{k-1}{2k}}}{\alpha^{\frac{15k+3}{4k}} \log^{\frac{k+1}{2k}} \rho}\right).$$

Let $B=\cup_{i=1}^t B_i$ and note that $|B|\leq pq =
\frac{n}{30\alpha}$. Repeating the argument of the proof of Lemma
\ref{large-minor} we obtain a subset $D$ such that the subgraph
$G''$ induce by $D$ has the following properties.

\begin{itemize}
\item
For every two disjoint subsets $X, |X|\leq
n/2$ and $Y,|Y|\leq \frac{1}{3\alpha}|X|$ of $G''$ we have that
$$e\big(X,V(G'')-(X\cup Y)\big) \geq \frac{r}{4\alpha}|X|.$$
\item
At least $m=\frac{t}{3\alpha}$ sets $B_j$ satisfy that $e(B_j,D) \geq \frac{r}{4\alpha}
|B_j|$.
\end{itemize}

Without loss of generality, we can assume that $B_1, \ldots, B_m$
satisfy: $e(B_j,D) \geq \frac{r}{4\alpha} |B_j|$. Let $U_j$ be the
set of neighbors of $B_j$ in $D$. Since our graph is $C_{2k}$-free
we have, by the result of Naor and Verstra\"{e}te \cite{NV}, that
$$e(B_j,U_j) \leq (2k-3)\Big((|B_j||U_j|)^{\frac{k+1}{2k}}+|B_j|+|U_j|\Big).$$
This inequality together with $k\geq 3$ and $|B_j|=q$ implies that
$|U_j| \geq (40k^2\alpha^{3/2}\log \rho)n/q$. Pick uniformly at
random with repetition $q/(40k^2\alpha^{3/2}\log
\rho)>r^{k/2-1}\log r$ vertices of $G''$ and denote this set by
$W$. For every index $1 \leq i \leq m$ the probability that $W$
does not intersect $U_i$ is at most
$\left(1-\frac{|U_i|}{|G''|}\right)^{|W|}\leq 1/e$. Therefore the
expected number of  sets $U_i$ that have a non-empty intersection
with $W$ is at least $(1-1/e)m<m/2$. Hence there is a choice of
$W$ that intersects at least $m/2\geq t/(8\alpha)$ sets $U_i$. By
Lemma \ref{diameter-C_2k}, $G''$ contains a connected subgraph on
$\leq (40k^2\alpha^{3/2}\log \rho)|W|\leq q$ vertices that
contains $W$. By adding extra vertices if necessary we can
construct a connected subset $B_{t+1}$ of size $q$ that contains
$W$ and hence is connected by an edge to at least $t/(8\alpha)$
sets $U_i, i \leq t$. This completes the proof of the lemma.
\hfill $\Box$

Substituting in the above lemma $\alpha=a^2$ and $\rho=e^{2a}$ we
obtain the following corollary.

\begin{coro}
\label{c51} Let $1 \leq a \leq \log r$, and let $G$ be a
$C_{2k}$-free graph of order $n$ such that $a^{26}r^k  \leq n \leq
e^{2a}r^k$ and for every two disjoint subsets $X, |X|\leq 0.7n$
and $Y,|Y|\leq \frac{1}{2a^2}|X|$ we have that
$e\big(X,V(G)-(X\cup Y)\big) \geq \frac{r}{4a^2}|X|$. Then $G$
contains a minor with average degree at least
$cr^{\frac{k+1}{2}}$, where $c$ is a constant independent of $r$
and $a$.
\end{coro}

\begin{lemma}
\label{l52}
Let $\rho\geq 3$ be a constant and let
$G$ be a $C_{2k}$-free graph on $\rho r^k$
vertices with average degree $r$.
Then $G$ contains a minor with average degree at least
$\Omega\big(r^{\frac{k+1}{2}}\big)$.
\end{lemma}

\noindent
{\bf Proof.}\,
Set $\alpha=8(k\log(32k)+\log \rho)$ and note that it is a constant independent of $r$.
If $G$ has the property that for every two disjoint subsets $X, |X|\leq
0.7n$ and $Y,|Y|\leq \frac{1}{2\alpha}|X|$ we have that
$$e\big(X,V(G)-(X\cup Y)\big) \geq \frac{r}{4\alpha}|X|,$$ then by
Lemma \ref{large-minor-C_2k} it contains a minor with average
degree $\Omega\left(r n^{\frac{k-1}{2k}}\right)$. Since by Lemma
\ref{C_2k-size} every $C_{2k}$-free graph with average degree
$\Omega(r)$ has at least $\Omega(r^k)$ vertices we are done.
Otherwise, there are two sets $X, Y$ as above for which
$e\big(X,V(G)-(X\cup Y)\big) < \frac{r}{4\alpha}|X|.$ Then, by
Lemma \ref{small-separation}, we have that the average degree of
the graph $G-X$ is at least $r$, or the average degree of the
subgraph induced by $X \cup Y$ is at least $r-\frac{r}{\alpha}$.
In the first case let $G_1=G-X$ and in the second let $G_1=G[X
\cup Y]$. Note that the number of vertices $n_1$ of the new graph
is strictly smaller than that of $G$. Moreover if the average
degree of $G_1$ is smaller than that of $G$ we know $n_1=|X\cup
Y|\leq 3n/4$. Continue this process until we either find a minor
with average degree at least
$\Omega\left(r^{\frac{k+1}{2}}\right)$, or we have at least
$\alpha/2$ steps at which the average degree of the new graph
decreases. In the second case, let $G'$ be the resulting graph and
$n'$ be the number of its vertices.

Since the degree decreased exactly $\alpha/2$ times we know that the average degree of $G'$ is at least
$r-(\alpha/2)\frac{r}{\alpha}\geq r/2$ and the number of its vertices satisfies
$$ n' \leq \left(\frac{3}{4}\right)^{\alpha/2} n < e^{-k\log(32k)-\log \rho}n=
\frac{n}{\rho (32k)^k} \leq \left(\frac{r}{32k}\right)^k.$$ As
$G'$ is $C_{2k}$-free, it contradicts the assertion of Lemma
\ref{C_2k-size}.  This shows that the second case is in fact
impossible and our process always outputs a minor of average
degree at least $\Omega\left(r^{\frac{k+1}{2}}\right)$. \hfill
$\Box$

\begin{lemma}
\label{smallgraphs-C_2k}
Let $G$ be a $C_{2k}$-free graph of with average degree $r$ and at most
$r^{k+2}$ vertices.
Then $G$ contains a minor with average degree at least $\Omega\left(r^{\frac{k+1}{2}}\right)$.
\end{lemma}

\noindent {\bf Proof.}\,  Let \{$a_i, i \geq 0\}$ be an increasing
sequence defined by  $a_0=65$ and $a_{i+1}=e^{a_i/13}$. Note that
$a^{26}_{i+1}=e^{2a_i}$ and let $\ell$ be the first index such
that $e^{2a_{\ell}} \geq r^{k+2}$. Then there is some $0 \leq i
\leq \ell$ the order $n$ of our graph $G$ satisfies $a_i^{26} r^k
\leq n < e^{2a_i}r^k$. If $G$ has the property that for every two
disjoint subsets $X, |X|\leq 0.7n$ and $Y,|Y|\leq
\frac{1}{2a_i^2}|X|$ we have that
$$e\big(X,V(G)-(X\cup Y)\big) \geq \frac{r}{4a_i^2}|X|,$$ then by
Corollary \ref{c51} it contains a minor with average degree
$\Omega\left(r^{\frac{k+1}{2}}\right)$ and we are done. Otherwise,
there are two sets $X, Y$ as above for which $e\big(X,V(G)-(X\cup
Y)\big) < \frac{r}{4a_i^2}|X|.$ Then, by Lemma
\ref{small-separation}, we have that the average degree of graph
$G-X$ is at least $r$, or the average degree of the subgraph
induced by $X \cup Y$ is at least $r-\frac{r}{a_i^2}$. In the
first case let $G_1=G-X$ and in the second let $G_1=G[X \cup Y]$.
Note that the number of vertices $n_1$ of new graph is strictly
smaller than that of $G$. Moreover if the average degree of $G_1$
is smaller than that of $G$ we know $n_1=|X\cup Y|\leq 3n/4$.
Continue this process until we either find a minor of with average
degree at least $\Omega\left(r^{\frac{k+1}{2}}\right)$ or we
arrive to a graph $G'$ with $n'$ vertices such that $n' \leq
a_0^{26}r^k$.

In the first case we are clearly done. In the second case we claim
that the average degree of $G'$ is still at least $r/2$. Note that
if at some stage the order of our graph $G_j$ satisfied $
a_i^{26}r^k \leq |V(G_j)| <e^{2a_i}r^k,$ then the average degree
of the new graph $G_{j+1}$ could decrease only by at most
$r/a_i^2$. In this case the order of $G_{j+1}$  drops as well so
that $|G_{j+1}| \leq 3|G_j|/4$. Since $(3/4)^4<e^{-1}$, we have
that this can happen only at most $8a_i$ times, before the order
of the remaining graph will become smaller than
$a_i^{26}r^k=e^{2a_{i-1}}r^k$. Since $a_{i+1}=e^{a_i/13}\geq
2a_i$,  we have that during all iterations the average degree of
the resulting graph can decrease by at most
$$ \sum_i8a_i \cdot \frac{r}{a_i^2}=r \sum_i \frac{8}{a_i} \leq \frac{16}{a_0} r <r/2.$$
Hence the final graph $G'$ has average degree at least $r/2$ and at most
$O(r^k)$ vertices. Therefore, by Lemma \ref{l52}, it contains a
minor with average degree $\Omega\left(r^{\frac{k+1}{2}}\right)$.
\hfill $\Box$

\vspace{0.25cm} \noindent {\bf Proof of Theorem
\ref{C_2k-minor}.}\, Let $G$ be a $C_{2k}$-free graph with average
degree $r$ and let $n$ be the number of vertices of $G$. By Lemma
\ref{smallgraphs-C_2k}, we can assume that $n > r^{k+2}$. Suppose
that $G$ contains a subset $X, |X| \leq 0.7n$, such that $|N(X)|
\leq \frac{|X|}{2\log^2 n}$. If the average degree of $G-X$ is at
least $r$, set $G_1=G-X$ and let $n_1$ be the number of vertices
in $G_1$. Otherwise, let $G_1$ be the subgraph induced by the set
$X\cup N(X)$. In the second case, by Lemma \ref{small-separation},
the average degree of $G_1$ is at least $r-\frac{r}{\log^2 n}$.
Note that in both cases we obtain a smaller graph. Moreover if the
average degree of $G_1$ is smaller than that of $G$ we know that
$n_1=|X\cup N(X)|\leq 3n/4$. Continue this process until we obtain
a subgraph $G'$ of $G$ on $n'$ vertices such that one of the
following holds. Either $n' \leq r^{k+2}$ or every subset $X$ of
$G'$ of size $|X| \leq 0.7n'$ has $|N(X)|  \geq \frac{|X|}{2\log^2
n'}$. Note that in the second case the graph $G'$ does not have a
separator of size $\frac{n'}{2\log^2 n'}$. Since $n'>r^{k+2}$, by
the result of Plotkin, Rao and Smith \cite{PRS}, $G'$ has a clique
minor of size
$$\Omega\left(\frac{n'/\log^2 n'}{\sqrt{n' \log n'}}\right)
=\Omega\left(\frac{\sqrt{n'}}{\log^{5/2} n'}\right) \geq
r^{\frac{k+2}{2}-o(1)} \gg r^{\frac{k+1}{2}}.$$ In the first case,
when $n' \leq r^{k+2}$ we claim that the average degree of $G'$ is
still at least $r/2$. Indeed, let $x_0=n, x_1, \ldots, x_\ell \geq
r^{k+2}$ be the sequence of orders of graphs that we had during
the process when the average degree decreased. Then we know that
$x_{i+1}\leq 3 x_i/4$ and the decrease in the average degree at
the corresponding step was at most $r/\log^2 x_i$. Let $y_i=\log
x_{\ell-i}$, then $y_0 \geq (k+2)\log r$ and $y_{i+1}\geq
y_i+\log(4/3) \geq y_i+1/4$. Therefore
$$ \sum_i \frac{1}{y^2_i} \leq \sum_i \frac{1}{(y_0+i/4)^2} \leq 16 \sum_{i=0}^{\infty}
\frac{1}{(4y_0+i)^2} \leq \frac{4}{y_0-1}\ll 1/2,$$ and we
conclude that the average degree of $G'$ is at least
$r\big(1-\sum_i 1/\log^2 x_i\big) \geq r/2$. Therefore we can find
in $G'$ a minor with average degree
$\Omega\left(r^{\frac{k+1}{2}}\right)$ using Lemma
\ref{smallgraphs-C_2k}. This completes the proof of the theorem.
\hfill $\Box$

\section{Concluding remarks}\label{s-final}

In  this paper we proved that if $G$ is an expander graph than it contains a large  
clique minor. Moreover our results on $H$-free graphs suggest that 
already local expansion may be sufficient to derive results of
this sort. This leads to the following 
general question which we think deserves further study. Let $G$ be a graph of order $n$ such that for 
every subset of vertices $X$ of size at most $s$ we have that 
$|N(X)| \geq t|X|$. Denote by
$f(s,t)$ the size of the largest clique minor which such graph must always contain.
What is the asymptotic behavior of this function? Note that we already know the behavior of $f$ in 
the two extremal cases when  $s=1$ and 
$s=\Theta(n/t)$. Indeed, if $s=1$ we just have that the minimum degree of $G$ is at least $t$ and therefore 
it contains a clique minor of order $\Omega\big(t/\sqrt{\log t}\big)$ by 
Kostochka-Thomason. In the second case we have by Theorem \ref{th1} that
our graph has clique minor of order
$\Omega\big(t\sqrt{s\log t}/\sqrt{\log(st)}\big)$.

One related and quite attractive question which remains unsettled
is the asymptotic behavior of the largest clique minor size in
sparse random graphs $G_{n,p}$. While for the case of constant
edge probability $p$, Bollob\'as, Catlin and Erd\H{o}s \cite{BCE}
showed this quantity to behave asymptotically as
$\Theta\big(n/\sqrt{\log n}\big)$, their method is apparently insufficient
to resolve the question for (much) smaller values of $p(n)$, and
in particular, for the the rather intriguing case $p=c/n$, $c>1$
is a constant, where a largest clique minor can be shown to be
with high probability between $c_1\sqrt{n/\log n}$ and
$c_2\sqrt{n}$.  

Another interesting direction of future study can be to find
sufficient conditions for ensuring a minor of a fixed graph $\Gamma$
(rather then just a clique $K_k$) in an expanding graph $G$. The
first step in this direction has been made by Myers and Thomason
\cite{MT05} who derived an analog of the Kostochka-Thomason result for a
general $\Gamma$.

Finally, it would be quite nice to obtain algorithmic analogs of
our main results (see, e.g. \cite{DHK} for a recent contribution
to algorithmic graph minor theory), providing efficient, deterministic 
algorithms for finding large minors, matching
our existential statements.

\end{document}